\documentclass[12pt,dvips]{amsart}
\usepackage{amsfonts, amssymb, latexsym, epsfig}
\usepackage{euler, amsfonts, amssymb, latexsym, epsfig, epic}

\newtheorem{Theorem}{Theorem}[section] 

\newtheorem{Proposition}[Theorem]{Proposition} 
\newtheorem{Lemma}[Theorem]{Lemma}

\newtheorem{Corollary}[Theorem]{Corollary}
\newtheorem{Problem}[Theorem]{Problem}
\newtheorem{Definition-Proposition}[Theorem]{Definition-Theorem}
\newtheorem{Main Conjecture}[Theorem]{Main Conjecture}

\newtheorem{Conjecture}[Theorem]{Conjecture}
\newtheorem{Definition}[Theorem]{Definition}

\theoremstyle{remark}
\newtheorem{Example}[Theorem]{Example}

\theoremstyle{plain}

\newcommand{\comment}[1]{$\star${\sf\textbf{#1}}$\star$}

\newcommand{\cellsize}{15}
\newlength{\cellsz} 
\newsavebox{\cell}
\sbox{\cell}{\begin{picture}(\cellsize,\cellsize)
\put(0,0){\line(1,0){\cellsize}}
\put(0,0){\line(0,1){\cellsize}}
\put(\cellsize,0){\line(0,1){\cellsize}}
\put(0,\cellsize){\line(1,0){\cellsize}}
\end{picture}}
\newcommand\cellify[1]{\def\thearg{#1}\def\nothing{}%
\ifx\thearg\nothing
\vrule width0pt height\cellsz depth0pt\else
\hbox to 0pt{\usebox{\cell} \hss}\fi%
\vbox to \cellsz{
\vss
\hbox to \cellsz{\hss$#1$\hss}
\vss}}
\newcommand\tableau[1]{\vtop{\let\\\cr
\baselineskip -16000pt \lineskiplimit 16000pt \lineskip 0pt
\ialign{&\cellify{##}\cr#1\crcr}}}
\newcommand{\kellsize}{30}
\newlength{\kellsz} 
\newsavebox{\kell}
\sbox{\kell}{\begin{picture}(\kellsize,\kellsize)
\put(0,0){\line(1,0){\kellsize}}
\put(0,0){\line(0,1){\kellsize}}
\put(\kellsize,0){\line(0,1){\kellsize}}
\put(0,\kellsize){\line(1,0){\kellsize}}
\end{picture}}
\newcommand\kellify[1]{\def\thearg{#1}\def\nothing{}%
\ifx\thearg\nothing
\vrule width0pt height\kellsz depth0pt\else
\hbox to 0pt{\usebox{\kell} \hss}\fi%
\vbox to \kellsz{
\vss
\hbox to \kellsz{\hss$#1$\hss}
\vss}}
\newcommand\ktableau[1]{\vtop{\let\\\cr
\baselineskip -16000pt \lineskiplimit 16000pt \lineskip 0pt
\ialign{&\kellify{##}\cr#1\crcr}}}

\newcommand{\sellsize}{36}
\newlength{\sellsz} 
\newsavebox{\sell}
\sbox{\sell}{\begin{picture}(\sellsize,20)
\put(0,0){\line(1,0){\sellsize}}
\put(0,0){\line(0,1){\sellsize}}
\put(\sellsize,0){\line(0,1){\sellsize}}
\put(0,\sellsize){\line(1,0){\sellsize}}
\end{picture}}
\newcommand\sellify[1]{\def\thearg{#1}\def\nothing{}%
\ifx\thearg\nothing
\vrule width0pt height\sellsz depth0pt\else
\hbox to 0pt{\usebox{\sell} \hss}\fi%
\vbox to \sellsz{
\vss
\hbox to \sellsz{\hss$#1$\hss}
\vss}}
\newcommand\stableau[1]{\vtop{\let\\\cr
\baselineskip -16000pt \lineskiplimit 16000pt \lineskip 0pt
\ialign{&\sellify{##}\cr#1\crcr}}}

\begin{document}
\pagestyle{plain}

\mbox{}
\title{Longest increasing subsequences, Plancherel-type measure\\ 
and the Hecke insertion algorithm}
\author{Hugh Thomas}
\address{Department of Mathematics and Statistics, University of New Brunswick, Fredericton, New Brunswick, E3B 5A3, Canada }
\email{hugh@math.unb.ca}

\author{Alexander Yong}
\address{Department of Mathematics, University of Minnesota, Minneapolis, MN 55455, USA}

\email{ayong@math.umn.edu}

\date{March 29, 2008}

\maketitle
\begin{abstract}
We define and study the \emph{Plancherel-Hecke probability 
measure} on Young diagrams; the \emph{Hecke algorithm}
of [Buch-Kresch-Shimozono-Tamvakis-Yong '06] 
is interpreted as a polynomial-time exact sampling algorithm for
this measure. Using the results of [Thomas-Yong '07] on \emph{jeu de taquin}
for \emph{increasing tableaux}, a symmetry
property of the Hecke algorithm is proved,
in terms of longest strictly increasing/decreasing subsequences
of words. This parallels 
classical theorems of [Schensted '61] and of 
[Knuth '70], respectively, on the \emph{Schensted}
and \emph{Robinson-Schensted-Knuth} algorithms. We investigate, and conjecture about, the limit 
typical shape of the measure, in analogy with work of [Vershik-Kerov '77],
[Logan-Shepp '77] and others on the ``longest increasing subsequence problem''
for permutations. We also include a related extension of 
[Aldous-Diaconis '99] on \emph{patience sorting}. Together, 
these results provide a new rationale 
for the study of increasing tableau combinatorics,
distinct from the original algebraic-geometric 
ones concerning  $K$-theoretic Schubert calculus.
\end{abstract}
\tableofcontents

\section{Introduction and main results}

\subsection{Overview} 
Let $W_{n,q}$ denote the set of words of length $n$ generated using the alphabet
$\{1,2,\ldots,q\}$. Let ${\tt LIS}(w)$ denote the {\bf length of the
longest strictly increasing subsequence} of $w=w_1 w_2\cdots w_n$, i.e., 
the largest $\ell$ with a subsequence
$i_1<i_2<\ldots<i_{\ell}$ such that 
$w_{i_1}<w_{i_2}<\ldots<w_{i_\ell}$. Similarly, we
consider the {\bf length of the
longest strictly decreasing subsequence} ${\tt LDS}(w)$ of~$w$. 
Our main goal is to 
introduce and study a discrete probability
measure on Young diagrams, in connection with the study of
the distributions of {\tt LIS} and {\tt LDS} on uniform random words.
An additional goal is to provide a novel motivation for the $K$-theoretic 
Schubert calculus combinatorics of \cite{BKSTY, Thomas.Yong:V}.

There are analogies with the study of {\tt LIS} and
{\tt LDS} in the {\bf permutation case}, i.e.,
when $w$ is chosen uniformly at random from the 
symmetric group $S_n$. The latter topic has attracted considerable attention; 
we refer the reader to the 
surveys \cite{Aldous.Diaconis, Stanley:inc} and the references therein.
In the permutation case, 
random Young diagrams are distributed according to the 
\emph{Plancherel measure} (on irreducible representations) of $S_n$.
This discrete probability measure is the push-forward of the
uniform distribution on $S_n$, under the \emph{Robinson-Schensted correspondence}. Schensted \cite{Schensted} established that this correspondence 
encodes ${\tt LIS}(w)$ and ${\tt LDS}(w)$ 
symmetrically in the shape $\lambda$ associated to 
$w$. In \cite{Vershik.Kerov, Logan.Shepp}, these ideas are applied
to determine the
asymptotics of the expectation of ${\tt LIS}$ over $S_n$ (solving the old
``longest increasing sequences problem''), via a study
of the ``limit typical shape'' under the Plancherel measure. 

As a continuation of this theme, we  
apply the {\tt Hecke} (insertion) algorithm of \linebreak
\cite{BKSTY} to define the Young diagram ${\tt Heckeshape}(w)$
for each $w\in W_{n,q}$; using this we 
define \emph{Plancherel-Hecke measure}.
Our belief that this measure should actually be worthy of analysis 
was initially guided by our theorem that
{\tt Hecke} symmetrically encodes ${\tt LIS}(w)$ and 
${\tt LDS}(w)$ for $w\in W_{n,q}$, a 
generalization of Schensted's theorem. During the course of
our investigation, we found that 
many other aspects of the 
Plancherel-Hecke measure (conjecturally) also resemble those of the 
Plancherel measure. This paper records these results, both
theoretical and computational, as a justification for further study.

Briefly, this is how the two aforementioned measures compare:
Let $q=\Theta(n^{\alpha})$. Sending $q,n\to\infty$, we conjecture that for
 $\alpha>\frac{1}{2}$, our measure is concentrated around 
the limit typical shape under Plancherel measure. This
\emph{Plancherel curve} plays an important role in
\cite{Vershik.Kerov, Logan.Shepp}. On the other hand, for
$\alpha<\frac{1}{2}$ we conjecture the measure is concentrated near the 
``staircase shape''. In particular, 
a ``phase transition'' is suggested at $\alpha=\frac{1}{2}$.
As we tune $\alpha$, a symmetric deformation of the Plancherel curve occurs.
In view of the above mentioned result on the {\tt Hecke} algorithm, this 
transition phenomenon is further evidenced by computations (with contributions
by O.~Zeitouni)
of the expectation 
of ${\tt LIS}$ and ${\tt LDS}$ as $\alpha$ varies; see Section~5 and the 
Appendix.

There have been earlier extensions of the permutation case
to $W_{n,q}$. The 
limit distribution of the length of the longest \emph{weakly} 
increasing/decreasing 
subsequence ({\tt LwIS}/{\tt LwDS}) on $W_{n,q}$ was found in work of 
\cite{Tracy.Widom}, following the breakthrough \cite{Baik.Deift.Jo}
on the limit distribution of ${\tt LIS}$ on $S_n$.  See also the more
recent work \cite{Houdre}.
However, analogous understanding of 
the distribution of {\tt LIS} and {\tt LDS} on $W_{n,q}$
appears to be less developed; 
see, e.g., \cite{Biane, Borodin.Olshanski, Tracy.Widom} for contributions.

As a point of comparison and contrast with our approach, previous
work on ${\tt LIS}$, ${\tt LDS}$ and $W_{n,q}$
utilizes the combinatorics of
the \emph{Robinson-Schensted-Knuth correspondence}, which
asymmetrically encodes ${\tt LwIS}$ and ${\tt LDS}$.
We offer an alternative viewpoint on the
relationship between Young diagrams and {\tt LIS}, {\tt LDS}. 
New questions and conjectures are raised, stemming from the Coxeter-theoretic 
viewpoint of \cite{BKSTY} (which in turn generalizes ideas of
\cite{Edelman.Greene}).

This text expresses our desire to point out a natural
link between the probabilistic combinatorics of {\tt LIS}, {\tt LDS} and
the combinatorial algebraic geometry of $K$-theoretic Schubert calculus. 
In particular, we apply and further develop the \emph{jeu de taquin}
for \emph{increasing tableaux} from \cite{Thomas.Yong:V}, thereby giving
another perspective on that work, distinct
from the original one. In summary, we believe that the availability of these
two disparate 
interpretations for \cite{BKSTY, Thomas.Yong:V}
provides something atypical to recommend $K$-theoretic 
tableau combinatorics, among the large array of 
interesting generalizations of the classical Young tableau and 
symmetric function theories known today.

\subsection{Plancherel-Hecke measure}
We identify a partition 
$\lambda=(\lambda_1\geq \lambda_2\geq \ldots \geq \lambda_k>0)$
with its Young diagram (in English notation); set 
$|\lambda|:=\sum_{i}\lambda_i$. Let 
${\mathbb Y}$ denote the set of all Young diagrams.
A filling of a shape $\lambda$ with a subset of the 
labels $\{1,2,\ldots,q\}$
is an {\bf increasing tableau} 
if it is strictly increasing in both rows and columns.
Let ${\tt INC}(\lambda,q)$ be the set of all increasing tableaux
of shape $\lambda$. 

\begin{figure}[t]
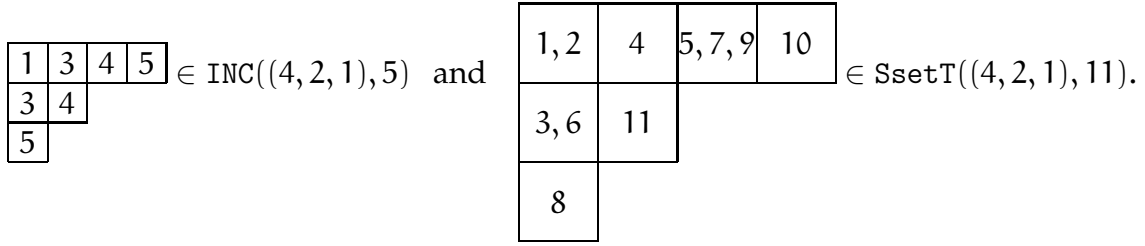

\[
\begin{tableau}
{{1} & {3} & {4} & {5}\\
{3} & {4}\\{5}}\end{tableau}\in {\tt INC}((4,2,1),5)
\mbox{\ \ \ and \ \ \ }
\begin{ktableau}
{{1,2}&{4}&{5,7, 9}&{10}\\
{3,6}&{11}\\
{8}}\in {\tt SsetT}((4,2,1), 11)
\end{ktableau}.\]
\caption{\label{fig:inc_set} An increasing tableau and a set-valued standard Young tableau}
\end{figure}

We also need {\bf set-valued tableaux} \cite{Buch:KLR}, which are
fillings of $\lambda$ assigning to each box
a nonempty subset of $\{1,2,\ldots,n\}$ such that
the largest entry of a box is smaller than the smallest entry in the
boxes directly to the right of it, and directly below it. We call a set-valued
tableau {\bf standard} if each label is used precisely once. Let
${\tt SsetT}(\lambda,n)$ denote the set of all standard set-valued tableaux.
See Figure~\ref{fig:inc_set}.

The {\bf Plancherel measure} on ${\mathbb Y}$ assigns to $\lambda$ the probability 
\[\frac{(f^{\lambda})^2}{n!}, \mbox{ \ where 
$f^{\lambda}:=e^{\lambda}(|\lambda|)$ is the
number of standard Young tableaux of shape $\lambda$.}\] 

Let 
$d^{\lambda}(q):=\#{\tt INC}(\lambda,q) \mbox{ \ and $e^{\lambda}(n)=\#{\tt SsetT}(\lambda,n)$}$. 

\begin{Definition}
\label{def:main}
The {\bf Plancherel-Hecke probability measure} $\mu_{n,q}$ on ${\mathbb Y}$
is defined by letting $\lambda_{n,q}$ be a random (non-uniform) Young 
diagram with distribution
\[Prob(\lambda_{n,q}=\lambda):=\frac{1}{q^n}d^{\lambda}(q)e^{\lambda}(n).\] 
\end{Definition}

\begin{Proposition}
\label{prop:central}
The Plancherel-Hecke measure is well-defined as a probability distribution;
i.e., the following identity holds: 
\begin{equation}
\label{eqn:ph2}
q^n=\sum_{\lambda}d^{\lambda}(q)e^{\lambda}(n),
\end{equation}
where 
\begin{equation}
\label{eqn:ph3}
|\lambda|\leq \min\left(n,{q+1\choose 2}\right) \mbox{\ \  and\ \  }
\lambda\subseteq (q, q-1,q-2,\ldots,3,2,1).
\end{equation}

There is an exact polynomial-time sampling algorithm 
\[{\tt Heckeshape}:W_{n,q}\to {\mathbb Y},\]
terminating in $O(nq^2)$ operations, 
that induces $\mu_{n,q}$ from the uniform distribution on $W_{n,q}$.
\end{Proposition}

The core technical result of this paper is a 
generalization of the aforementioned theorem of
Schensted~\cite{Schensted}:

\begin{Theorem}
\label{thm:hard}
{\tt Heckeshape} simultaneously and symmetrically encodes 
${\tt LIS}(w)$ and ${\tt LDS}(w)$ as the
size of the first row and column
of ${\tt Heckeshape}(w)$, respectively. 
\end{Theorem}

Theorem~\ref{thm:hard} is obtained by establishing another new result, 
connecting ${\tt Heckeshape}$ to the ``$K$-infusion'' operation defined in
\cite{Thomas.Yong:V}. This latter result is an analogue of the classical
fact that connects the Robinson-Schensted correspondence to the (ordinary)
\emph{jeu de taquin} rectification procedure.

We prove of Proposition~\ref{prop:central} in Section~2,
after recalling the {\tt Hecke} algorithm of 
Buch-Kresch-Shimozono-Tamvakis-Yong 
\cite{BKSTY} (originally
constructed to study degeneracy loci of vector bundles).
${\tt Heckeshape}(w)$ is the Young diagram associated to $w$ under
${\tt Hecke}$. The proof of
Theorem~\ref{thm:hard} is given in
Section~4. 


\begin{Example}
We illustrate the identity (\ref{eqn:ph2}) for 
$n=4$ and $q=3$. There are nine partitions $\lambda$ satisfying 
(\ref{eqn:ph3}). These are
\[(1), \ (2),\  (1,1),\  (2,1),\ (3),\ (1,1,1),\ (3,1),\ (2,1,1),\  (2,2).\]
Then (\ref{eqn:ph2}) reads
\[81=3^4=3\cdot 1+3\cdot 3+3\cdot 3+5\cdot 8+1\cdot 3+1\cdot 3+2\cdot 3+2\cdot 3+1\cdot 2,\]
where the products on the righthand side of the equality are listed in
order corresponding to the above partitions. Thus, the ``typical shape''
is $(2,1)$, possessing nearly $50\%$ of the distribution. The remainder of 
the above results will be illustrated in Section~2, 
after we define
{\tt Heckeshape}.
\end{Example}

Theorem~\ref{thm:hard} has some immediate consequences, familiar from the
permutation case.

\begin{Corollary}
\label{cor:expectformula}
Under the uniform measure on $W_{n,q}$, we have
\begin{equation}
\label{eqn:expectformula}
E({\tt LIS})=\frac{1}{q^n}\sum_{\lambda}\lambda_1 d^{\lambda}(q)e^{\lambda}(n).
\end{equation}
In addition,
\[Prob({\tt LIS}=\ell)=\frac{1}{q^n}\sum_{\lambda, \lambda_1=\ell}d^{\lambda}(q)e^{\lambda}(n).\]
\end{Corollary}

Two other consequences will be given in Section~3. 
One gives a ``Coxeter-theoretic'' generalization of 
the widely known Erd\H{o}s-Szekeres theorem \cite{Erdos}. 
Another expands upon the discussion of 
patience sorting given in \cite{Aldous.Diaconis}.

\subsection{Remarks on Proposition~\ref{prop:central} and 
Theorem~\ref{thm:hard}}

In our experiments, {\tt Heckeshape} was reasonably efficient as a 
sampling algorithm.\footnote{Software available at the
authors' websites.} For example, when $n\leq 10,000$ sampling one Young
diagram takes on the order of seconds to minutes on current technology. 
For larger $n$, we could sample one
Young diagram when $n=50,000$ in several hours on the same 
technology. A sample when $n=100,000$ took about one and a half days.
The memory demands were modest. In view of the apparent ``concentration'' 
suggested below, one sample was enough to be of interest for our purposes, when
$n$ is large. 

There are classical antecedents of 
Theorem~\ref{thm:hard}. As stated earlier, 
Schensted \cite{Schensted} proved the 
analogous conclusion about the shape coming from the
Robinson-Schensted correspondence for 
a permutation $w$. In contrast, Knuth \cite{Knuth}
proved that the first row of the Robinson-Schensted-Knuth algorithm
({\tt RSK}) encodes the length of the 
longest weakly increasing subsequence (${\tt LwIS}$) 
of a word~$w$. 

What is perhaps less well-known is that {\tt RSK}
also encodes ${\tt LDS}(w)$ as the length of the first
column of the shape it associates to $w$. 
However, unlike {\tt Heckeshape},
it is asymmetric:
${\tt LIS}(w)$ is not encoded by the length of the first row
(as ${\tt LIS}(w)\neq {\tt LwIS}(w)$ in general). Thus, the symmetry of
{\tt Hecke(shape)} makes it natural to analyze, 
as it seems desirable to simultaneously capture the
statistics of ${\tt LIS}$ and ${\tt LDS}$.

However, we do not have handy formulas for $Prob(\lambda)$,
$d^{\lambda}(q)$ or $e^{\lambda}(n)$, such as the
\emph{hook-length formula} for $f^{\lambda}$.
(In the permutation case, the hook-length formula plays a crucial
role, see \cite{Vershik.Kerov, Logan.Shepp}.) 
In small examples large prime factors appear, 
showing that such a formula is unlikely. For instance, 
$d^{(4,2,1)}(7)=1337=7\cdot 191 \mbox{\  and\  } e^{(4,2,1)}(8)
=452=2^2\cdot 113$. This issue is closely related to the open question
of finding ``good'' 
determinantal expressions for \emph{Grothendieck polynomials}~\cite{LS:Hopf}.
That being said, 
special cases exhibit connections to work of Stanley \cite{Stanley:poly}
on polygonal dissections, and of \cite{Fomin.Greene} on generalized 
Littlewood-Richardson rules (further discussion may appear elsewhere). 
Thus, the enumerative combinatorics of these numbers might of interest in their
own right.
 

It is not difficult to give recursions to calculate $d^{\lambda}(q)$ 
and $e^{\lambda}(n)$ that are useful
in moderately large cases. Are there efficient (possibly randomized or approximate) counting algorithms?

Objectively, the lack of simple formulas to compute 
$Prob(\lambda)$ makes it trickier to apply standard approaches
directly; this is an admitted defect of our setting. 
Nevertheless, we believe the framework of problems described here 
is tractable. 
In addition to the results below, in the Appendix one gains useful and 
nontrivial information about the Plancherel-Hecke measure by 
exploiting the related work of~\cite{Biane}. In this way, the techniques 
of \cite{Logan.Shepp, Vershik.Kerov} can be applied to the present
context. 

\subsection{Analysis of $\mu_{n,q}$ and the limit typical shape}
We organize our analysis by first setting
\begin{equation}
\label{eqn:setq}
q=f(n)\in \Theta(n^{\alpha}), \mbox{ \ \  \ where $0< \alpha\leq 1$,}
\end{equation} 
and considering the limit behavior of $\mu_{n,q}$ 
when $q\to \infty$, as $n\to \infty$. 
(The case $\alpha=0$ is trivial.) As is explained below,
we conjecture that there is a 
{\bf critical value} of $\alpha$, denoted
$\alpha_{\rm critical}:=\frac{1}{2}$: the behavior of
$\mu_{n,q}$ is qualitatively different in the intervals
$\alpha\in(\alpha_{\rm critical},1]$ and $\alpha\in (0,\alpha_{\rm critical})$.
At $\alpha=\alpha_{\rm critical}=\frac{1}{2}$, further refinement of the analysis 
is needed, as we transition from one state to the other.

In the permutation case, to study the Plancherel measure, it is 
useful to consider the most likely, or ``typical'' shape. There, 
three facts are true. First, in the large limit (and after
rescaling), a well-defined typical shape exists.
Second, the expectation of 
the {\tt LIS} and {\tt LDS} of a large random permutation 
is encoded respectively in the length of the first row and column of the 
limit shape. 
Third, 
the Plancherel measure is concentrated near the typical shape. 

We conjecture that analogues of
all three of the aforementioned features also 
hold for the Plancherel-Hecke measure.

To be more precise, let the {\bf typical shape}
$\Lambda_{n,q}$ be the shape
$\lambda$ (contained in $(q,q-1,\ldots,2,1)$)
maximizing $Prob(\lambda)$. This Young diagram $\Lambda_{n,q}$ can
be interpreted as a step-function. It can furthermore be associated
to a piecewise linear approximation 
$f_{n,q}:[0,\infty)\to {\mathbb R}_{\geq 0}$. Finally, rescale by
\[{\hat f}_{n,q}(x):=\left\{\begin{array}{cc}
\frac{1}{2\sqrt n}f_{n,q}(x\cdot{2\sqrt n}) & \mbox{if $\alpha\geq \alpha_{\rm critical}=\frac{1}{2}$}\\  
\frac{1}{q}f_{n,q}(x\cdot {q}) & \mbox{otherwise.}
\end{array}\right.\]
 
\begin{Conjecture}
\label{conj:thebigone}
\begin{itemize}
\item [(I)] For
any $0\leq \alpha\leq 1$, there is a unique continuous function
\[\Lambda\in C([0,\infty)\to {\mathbb R}_{\geq 0})\] 
such that for any $\epsilon>0$, 
\[Prob\left(\sup_{x\in {\mathbb R}_{\geq 0}} |{\hat f}_{n,q}-\Lambda|>\epsilon\right)\to 0\]
as $n\to \infty$. We call this $\Lambda$ the {\bf limit typical
shape}.
\item[(II)] A ``phase transition'' occurs at 
$\alpha_{\rm critical}=\frac{1}{2}$:

For $0< \alpha<\alpha_{\rm critical}=\frac{1}{2}$, $\Lambda$
is the line 
\begin{equation}
\label{eqn:line}
y=1-x  \mbox{ \ \ \ for $0\leq x\leq 1$.}
\end{equation}

For $\alpha_{\rm critical}=\frac{1}{2}< \alpha\leq 1$, 
$\Lambda$ is the {\bf Plancherel curve}, which is parametrically given by
\begin{equation}
\label{eqn:Vershikcurve}
x=y+\cos\theta, \  \ \ y=\frac{1}{\pi}(\sin\theta - \theta\cos\theta), \mbox{ \ \ \ for $0\leq \theta\leq\pi$.}
\end{equation}
(The curves defined by (\ref{eqn:line}) and (\ref{eqn:Vershikcurve})
are declared to be identically $0$ for $1\leq x<\infty$.)
\item[(III)] For $\alpha=\alpha_{\rm critical}=\frac{1}{2}$: there is a 
constant $C> 0$ such that if 
\[q=kn^{\frac{1}{2}}+\mbox{lower order terms},\] 
then if $k< C$, $\Lambda$ is given by (\ref{eqn:line}). 
Otherwise, $\Lambda$ is given by a
deformation of 
(\ref{eqn:Vershikcurve}) which is symmetric across the line $y=x$. 
In either case, the $x$ and $y$ intercepts are at 
\[0\leq\beta(k)\leq 1\] 
where
\[E({\tt LIS})\approx\beta(k){2\sqrt n}\] 
and explicitly,
\begin{equation}\nonumber
\beta(k)=\left\{
\begin{array}{cc}
\frac{k}{2} & \mbox{if $0<k\leq 1$}\\
\frac{2-k^{-1}}{2} & \mbox{if $k>1$}
\end{array}\right..
\end{equation}

\end{itemize}
\end{Conjecture}

We have reasonable support
for the cases (I) and (II) of Conjecture~\ref{conj:thebigone}. 
\emph{Heuristically}, part (II)
of the conjecture says that when $\alpha$ is large, and thus
$q$ is ``close'' to $n$, a random word is ``close to'' being a 
random permutation. For a permutation, Schensted and {\tt Hecke}
behave the same. Hence the Plancherel and Plancherel-Hecke measures ought to be
 maximized on the same shape. When $\alpha$ is small, the limit
typical shape is a rescaling of the ``staircase shape'' which plays a
distinguished role in the Edelman-Greene algorithm \cite{Edelman.Greene}
and the {\tt Hecke} algorithm; see Theorem~\ref{thm:desc} and its proof.

Conjecture~\ref{conj:thebigone}(III) is
more speculative, since we did not have as much computational
evidence for the shape of $\Lambda$. There appears to be 
a continuous ``flattening''
of the Plancherel curve to a line as we tune $k$ from $\infty$ to $0$.
Our data was insufficient to rule out the possibility
that $\Lambda$ is simply a rescaling of the Plancherel curve by a factor of $\beta(k)$. 

\begin{Problem}
Explicitly describe the deformation of $\Lambda$ when 
$\alpha=\alpha_{\rm critical}=\frac{1}{2}$, as $k$ varies.
\end{Problem}
Our best estimate 
is that $\frac{1}{2}\leq C\leq 1$ (probably just $C=1$). However, the values
of $\beta(k)$ for $k$ relatively small can be experimentally 
estimated. The table
below was based on Monte Carlo estimates of $E({\tt LIS})$ for 
$n=50,000,\  100,000,\  200,000$ and $300,000$ and the estimates
were stable throughout this range. They closely 
agree with the conjecture for $\beta(k)$ given above.

\begin{table}[h]
\begin{tabular}{|l||l|l|l|l|l|}
\hline
$k$ & 0.5 & 1 & 2 & 4 & 10\\ \hline
$\beta(k)$ estimate & 0.25 & 0.50 & 0.74 & 0.86 & 0.94 \\ \hline 
\end{tabular}
\caption{Estimates of $\beta(k)$ 
for the $\alpha=\alpha_{\rm critical}=\frac{1}{2}$ case}
\end{table}


Notice that 
since $Prob(\lambda)=Prob(\lambda')$ where $\lambda'$ is the conjugate shape
of $\lambda$ we know that
if $\Lambda$ exists and is unique, then $\Lambda$ is symmetric.
This is consistent with the limit curves we predict.


We prove in Section~5 that: 
\begin{Theorem}
\label{thm:desc}
Conjecture~\ref{conj:thebigone} is
true for $0\leq \alpha<\frac{1}{3}$. 
More precisely, in this range, a random shape $\lambda$ under
$\mu_{n,q}$ satisfies $\lambda_i=q-i+1$ almost surely, as $n,q\to \infty$.
\end{Theorem}

The proof of this theorem depends
on the analysis of a certain random walk on the symmetric group. 
Our analysis is not sharp enough to 
extend to the range $\frac{1}{3}\leq\alpha\leq \alpha_{\rm critical}=\frac{1}{2}$,
although a refinement might be possible.

Empirically, one finds that the first row and column of 
$\Lambda_{n,q}$ are approximations of $E({\tt LIS})$ and 
$E({\tt LDS})$ that improve as $n\to\infty$. 
Therefore, it makes sense to study the asymptotics of
$E({\tt LIS})$ and $E({\tt LDS})$ as a means to understand the characteristics
of $\Lambda_{n,q}$. From this point of view, the 
following result supports the phase
transition phenomena asserted in 
Conjecture~\ref{conj:thebigone}. 

\begin{Theorem}[With O.~Zeitouni]
\label{thm:main}
If $0\leq \alpha< \alpha_{\rm critical}=\frac{1}{2}$ then 
$\lim_{n\to \infty} E({\tt LIS})=q$, 
whereas if $\alpha_{\rm critical}=\frac{1}{2}<\alpha\leq 1$ then
$E({\tt LIS})\approx 2\sqrt n$. The same statements hold when
$E({\tt LDS})$ replaces $E({\tt LIS})$.
\end{Theorem}

The proof for $\alpha<\alpha_{\rm critical}=\frac{1}{2}$ is a
variation on the approach we use to prove Theorem~\ref{thm:desc}.
We also conjectured the answer for
$\alpha>\alpha_{\rm critical}=\frac{1}{2}$; after showing O.~Zeitouni our guess during an early
stage in the project, he communicated
to us a proof, and kindly allowed us to reproduce his argument here.

In private communication, E.~Rains offered a proof that
$E({\tt LIS})=q$ in the $\alpha=\alpha_{\rm critical}=\frac{1}{2}$ and $k\leq 1$ case.
Afterward, in the appendix for this paper,
O.~Zeitouni and the second author present a simple proof that 
$E({\tt LIS})\approx \beta(k)2\sqrt n$, for all $k$, in the
$\alpha=\alpha_{\rm critical}=\frac{1}{2}$ case, thereby closing the gap in 
Theorem~\ref{thm:main}. The proof builds on work of \cite{Biane} (see further
discussion in Section 1.5).
These results further support
the belief that $C=1$ in Conjecture~\ref{conj:thebigone}(III).


In addition, we have the following 
conjecture about the fluctuation of ${\tt LIS}$ and 
${\tt LDS}$.
\begin{Conjecture} 
\label{conj:lambda1equals}
Let $\sigma({\tt LIS})$ denote the standard deviation 
of ${\tt LIS}$. 
For $0<\alpha<\alpha_{\rm critical}=\frac{1}{2}$ then 
\[\lim_{n\to\infty} \sigma({\tt LIS})=0,\] 
whereas if $\alpha_{\rm critical}=\frac{1}{2}<\alpha\leq 1$ then
\[\lim_{n\to\infty} \sigma({\tt LIS})=O(n^{\frac{1}{6}}).\]
The same statements hold for ${\tt LDS}$.
\end{Conjecture}

Note that Theorem~\ref{thm:desc} implies 
Conjecture~\ref{conj:lambda1equals} holds 
for $0\leq \alpha<\frac{1}{3}$. Tables~\ref{table:50000MC} and~\ref{table:100000MC} give numerical evidence for Conjecture~\ref{conj:lambda1equals} and
are consistent with Theorem~\ref{thm:main}.

\begin{table}[h]
\caption{\label{table:50000MC} $n=50,000$ with $1,000$ Monte Carlo trials}
\begin{tabular}{|l||l|l|l|l|l|l|}
\hline
estimate $\backslash \alpha$ & 0.45 & 0.50 & 0.55 & 0.60 & 0.75 & 1.00\\ \hline
$E({\tt LIS})$ & 130.00 & 222.50 & 311.06 & 368.38 & 422.48 & 436.36 \\ \hline
$\sigma$ & 0.00 & 0.63 & 2.86 & 4.01 & 5.07 & 5.09\\ \hline
\end{tabular}

\bigskip

\caption{\label{table:100000MC} $n=100,000$ with $500$ Monte Carlo trials}
\begin{tabular}{|l||l|l|l|l|l|l|}
\hline
estimate $\backslash \alpha$ & 0.45 & 0.50 & 0.55 & 0.60 & 0.75 & 1.00\\ \hline
$E({\tt LIS})$ & 177.00 & 315.43 & 448.2 & 523.63 & 603.78 & 619.64 \\ \hline
$\sigma$ & 0.00 & 0.67 & 3.52 & 4.62 & 5.90  & 6.29 \\ \hline
\end{tabular}
\end{table}

\begin{figure}[h]
  \centering
  \epsfig{file=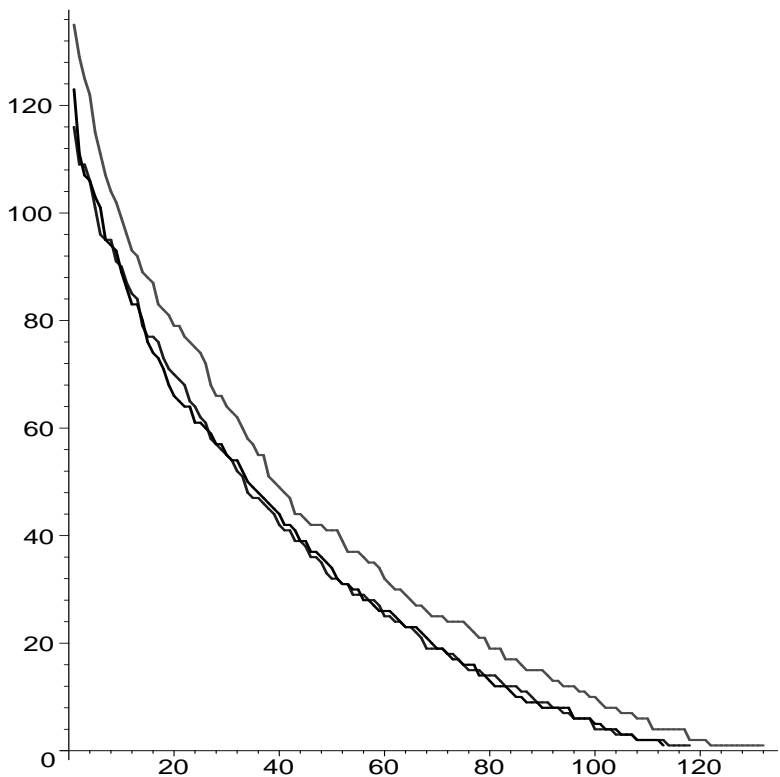,height=4in, width=4in}
\caption{\label{fig:converge} Two samples at $n=5000, \alpha=\frac{2}{3}$ compared with an empirical approximation of the Plancherel curve; conjecturally as $n\to\infty$, the sample curves converge to one another}
\end{figure}

The bulk of 
$\mu_{n,q}$ appears ``concentrated'' near $\Lambda_{n,q}$, i.e., 
the probability of sampling a random shape differing, in the sup-norm, 
from $\Lambda$ after rescaling, by some fixed $\epsilon>0$, goes to $0$ as 
$n,q\to \infty$. See Figure~\ref{fig:converge}: already at $n=5,000$ we see that two random samples are visibly
``close'' to one another, and are similar in shape 
to the third curve which is an
approximation of the Plancherel curve. 
By $n=100,000$ the curves appear undeniably to be rescalings of one 
another, with a rescaling factor of about $0.95$ in our experiments.
Naturally, as $\alpha$ gets larger, the empirical 
convergence of the curves occurs faster.

\subsection{Further comparisons with the literature}

As mentioned earlier, the limit distribution of ${\tt LIS}$
on permutations, and that of ${\tt LwIS}$ on words is well understood.

The study of {\tt LIS} on $W_{n,q}$ was considered, 
e.g., in~\cite{Tracy.Widom}. In addition, the study of the distribution of
{\tt LIS}, in the critical case $\alpha_{\rm critical}=\frac{1}{2}$ is implicit
in \cite{Biane, Borodin.Olshanski}. In \cite{Biane}, an alternative measure
on Young diagrams is studied: 
{\bf Schur-Weyl duality} implies that one has the decomposition
\[({\mathbb C}^{q})^{\otimes n}\cong\bigoplus_{\lambda}S^{\lambda}\otimes V_{\lambda},\]
where here $S^{\lambda}$ is the $S_n$ irreducible \emph{Specht module} and 
$V_{\lambda}$ is the $GL_{q}({\mathbb C})$ irreducible \emph{Schur module}. 
Now taking dimensions one defines a probability measure that assigns 
to $\lambda$ the likelihood 
$({\rm dim} \ S^{\lambda} \cdot {\rm dim} \ V_{\lambda})/q^n$.
Biane explicitly determines the rescaled limit typical shape in this context. 
Combinatorially, Biane's measure arises from the {\tt RSK} algorithm.
Since we know that ${\tt RSK}(w)$ encodes the ${\tt LIS}(w)$ in the 
first column (by reading $w$ backwards), 
one expects, by analogy with \cite{Logan.Shepp, Vershik.Kerov} that a 
certain rescaling of the first column of Biane's limit shape 
is the $\beta(k)$ of 
Conjecture~\ref{conj:thebigone}.
However, to justify this conclusion rigorously one needs more work. 
Further, the fluctuations around Biane's curve have been studied, in the 
$k=1$ case, by Borodin-Olshanski \cite{Borodin.Olshanski}.

{\tt Hecke} was originally
developed in \cite{BKSTY} as a generalization of the \emph{Edelman-Greene
correspondence} which bijects Coxeter reduced words in the symmetric group
to pairs of tableaux \cite{Edelman.Greene}. 
Our proof of Theorem~\ref{thm:hard} implies
that this algorithm encodes the {\tt LIS} of such words, although
the study of {\tt LIS} of reduced words appears unmotivated. On the other
hand, the Coxeter-theoretic viewpoint on words 
will be useful in our analysis of {\tt LIS}, {\tt LDS} and $\Lambda$.

\subsection{Summary and organization}
In Section~2 we recall the {\tt Hecke} algorithm and give an additional 
example of the results of Section~1.2. We then prove Proposition~\ref{prop:central}.
In Section~3, we include two 
consequences of Theorem~\ref{thm:hard}. We split our remaining 
proofs according to 
the main flavor of technique used: in Section~4 we explain the increasing
tableau theory we need from \cite{Thomas.Yong:V} 
and prove Theorem~\ref{thm:hard}. In Section~5, we utilize 
probabilistic-combinatorial techniques, combined with our main results,
to prove Theorems~\ref{thm:desc} and~\ref{thm:main}.

\section{The {\tt Hecke} algorithm}

\subsection{The $0$-Hecke monoid}
We need to recall some notions used in \cite{BKSTY}.
The {\bf $0$-Hecke monoid} ${\mathcal H}_{0,q}$ is the
quotient of the free monoid of all finite words in the alphabet
$\{1,2,\ldots,q\}$ by the relations
\begin{align}
\label{E:idem}
i\,i &\equiv i &\quad&\text{for all $i$} \\
\label{E:braid}
i\,j\,i &\equiv j\,i\,j &\quad&\text{for all $i,j$} \\
\label{E:comm}
i\,j &\equiv j\,i &\quad&\text{for $|i-j|\ge 2$.}
\end{align}

There is a bijection between ${\mathcal H}_{0,q}$ and the
symmetric group $S_{q+1}$. 
Given any word $a\in {\mathcal H}_{0,q}$
there is a unique permutation $\pi\in S_{q+1}$ such that $a\equiv b$
for any reduced word $b$ of $\pi$; see, e.g., the textbook
\cite{Bjorner.Brenti} for basic Coxeter theory for the symmetric group.  
In this case, we write $W(a)=\pi$ and say that $a$ is a {\bf Hecke word} for
$\pi$.  Indeed, the reduced words for $\pi$ are precisely the
Hecke words for $\pi$ that are of the minimum length $\ell(\pi)$, the 
{\bf Coxeter length} of $\pi$ (after identifying the label $i$ with the
simple reflection $s_i=(i \ \  i+1)$). Given an additional
permutation $\rho$ with Hecke word $b$, the {\bf Hecke product} of
$\pi$ and $\rho$ is defined as the permutation $\pi \cdot \rho =W(ab)$.

The (row reading) word of a tableau $T$, denoted ${\tt word}(T)$, is
obtained by reading the rows of the tableau from left to right,
starting with the bottom row, followed by the row above it, etc.  
We also define $W(T) := W({\tt word}(T))$.
So for example, if $T$ is the increasing tableau of Figure~\ref{fig:inc_set}, 
then ${\tt word}(T)=5\ 3 \ 4 \ 1\ 3\ 4 \ 5$.

The {\tt Hecke} algorithm defined in \cite{BKSTY} 
identifies pairs $(w,i)$ of words
\[w=w_1 w_2 \cdots w_n, \ \ \ \ i=i_1 i_2 \cdots i_n,\]
where $w$ is a {\bf Hecke word} and $i$
satisfies 
\[i_1\leq i_2\leq \ldots \leq i_n  \mbox{ \ \ and $i_j<i_{j+1}$ whenever
$w_j\leq w_{j+1}$},\] 
with pairs of tableaux $(P,Q)$ of the same shape, where $P$ is an increasing
tableau such that ${\tt word}(P)\equiv w$ and the {\bf content} (i.e., multiset
of labels) of $Q$ matches the content of $i$. 
We refer to the $P$-tableau as the {\bf insertion tableau} and the $Q$-tableau
as the {\bf recording tableau}. (We point out that the ``column'' 
convention in \cite{BKSTY} differs slightly from the ``row'' one used here.)

\subsection{Description of {\tt Hecke} and {\tt Heckeshape}} 
The following description of {\tt Hecke} was originally given in \cite{BKSTY}:

\noindent
\emph{Description of {\tt Hecke}:} In this algorithm, one 
inserts an integer $x$ into 
an increasing tableau $T$. We denote this by $T\leftarrow x$.
The output is a triple $(U,c,\alpha)$ 
where $U$ is a modification of $T$ (possibly $T=U$), $c$ is a corner
of $U$ and $\alpha\in\{0,1\}$ is a parameter. Initially, we attempt to 
insert $x$ into the first row of $T$, and an output integer is possibly
created which is inserted into the next row and so on, until no output
integer is created.  We refer to this final insertion as the {\bf 
terminating step}, and the previous insertions as {\bf bumping steps}.

Suppose $R$ is a row that we are attempting to insert $x$ into. If $x$
is larger than or equal to all the entries of $R$, then no output integer
is generated and the algorithm terminates: if adjoining $x$ to the end
of $R$ results in an increasing tableau $U$, then set $\alpha=1$ and $c$
to be the new corner added. Otherwise end with the present $U$, without
modification; $\alpha=0$ and $c$ is the corner that is at the end of the
column containing the rightmost box of $R$.
On the other hand, if $R$ contains boxes strictly larger than $x$, let
$y$ be the smallest such box. If replacing $y$ with $x$ results in an 
increasing tableau, then do so. In either case, $y$ is the output integer to
be inserted into the next row.

Inserting a word $w$ using this algorithm terminates with an increasing
tableau 
\[P=((((\emptyset\leftarrow w_1)\leftarrow w_2)\leftarrow \cdots \ )\leftarrow w_n).\]
The $Q$ tableau is obtained by
placing each $i_j$ in the $c$-corner resulting from the insertion of $w_j$.\qed

We also have the following reverse insertion algorithm ${\tt Hecke}^{-1}$.

\noindent
\emph{Description of ${\tt Hecke}^{-1}$:}
Let $Z$ be an increasing tableau, $c$ a corner of $Z$, and $\alpha \in
\{0,1\}$.  Reverse insertion applied to the triple
$(Z,c,\alpha)$ produces a pair $(Y,x)$ of an increasing tableau $Y$
and a positive integer $x$ as follows.  Let $y$ be the integer in the
cell $c$ of $Z$.  If $\alpha=1$, remove $y$.  In any case, reverse
insert $y$ into the row above the corner $c$.

Whenever a value $y$ is reverse inserted into a row $R$, let $x$ be
the largest entry of $R$ such that $x < y$.  If replacing $y$ with $x$
results in an increasing tableau, then this is done.  In any case, the
integer $x$ is passed up.  If $R$ is not the top
row, this means that $x$ is reverse inserted into the row above
of $R$; otherwise $x$ becomes the final output value, along with the
modified tableau.

We now complete the description of ${\tt Hecke}$.  Locate the bottom-most
corner with the largest label, in the $Q$ tableau, and remove the label.  
If it was the only
entry in its corner, remove the corner, set $\alpha=1$.  Otherwise set
$\alpha=0$.  Set $c$ to be this corner.  Then reverse insert $(P,c,\alpha)$.
Repeat until all the entries of $Q$ (and $P$) have been removed.  
\qed

{\tt Hecke} is a generalization of the Robinson-Schensted correspondence
in the sense that it agrees with that correspondence whenever $w$ is
a permutation in $S_n$. In that case the $P$ and $Q$ tableaux are both
standard Young tableaux.

In this paper, we are only concerned with the case $i_j=j$. 
Therefore, we also set ${\tt Hecke}(w):={\tt Hecke}(w,123\cdots n)$
and define 
\[{\tt Heckeshape}:W_{n,q}\to {\mathbb Y}\]
by setting ${\tt Heckeshape}(w)$ to be the common shape of $P$ and $Q$
under ${\tt Hecke}(w)$. (An alternative description of this map is given in
Theorem~\ref{thm:hecke_is_rect} in Section~4.)

\begin{Example}
Let $w= 5 \ 4 \ 1 \ 3 \ 4 \ 2 \ 5 \ 1 \ 2 \ 1 \ 4 \ 2 \ 4\in W_{13,5}$. 
Then the reader can check that
{\tt Hecke} produces the following steps:
\[\tableau{{5}},\tableau{{1}}\mapsto
\tableau{{4}\\{5}},\tableau{{1}\\{2}}\mapsto
\tableau{{1}\\{4}\\{5}},\tableau{{1}\\{2}\\{3}}\mapsto
\tableau{{1}&{3}\\{4}\\{5}},\tableau{{1}&{4}\\{2}\\{3}}\mapsto
\tableau{{1}&{3}&{4}\\{4}\\{5}},\tableau{{1}&{4}&{5}\\{2}\\{3}}
\mapsto\tableau{{1}&{2}&{4}\\{3}\\{4}\\{5}},\tableau{{1}&{4}&{5}\\{2}\\{3}\\{6}}\]
\[\mapsto
\tableau{{1}&{2}&{4}&{5}\\{3}\\{4}\\{5}},\tableau{{1}&{4}&{5}&{7}\\{2}\\{3}\\{6}}\mapsto
\tableau{{1}&{2}&{4}&{5}\\{2}\\{3}\\{4}\\{5}},\tableau{{1}&{4}&{5}&{7}\\{2}\\{3}\\{6}\\{8}}
\mapsto
\tableau{{1}&{2}&{4}&{5}\\{2}&{4}\\{3}\\{4}\\{5}},\tableau{{1}&{4}&{5}&{7}\\{2}&{9}\\{3}\\{6}\\{8}}\]
\[\mapsto
\tableau{{1}&{2}&{4}&{5}\\{2}&{4}\\{3}\\{4}\\{5}},\ktableau{{1}&{4}&{5}&{7}\\{2}&{9}\\{3}\\{6}\\{8,10}}\mapsto\cdots\mapsto\]
\[\tableau{{1}&{2}&{4}&{5}\\{2}&{4}&{5}\\{3}&5\\{4}\\{5}},
\stableau{{1}&{4}&{5}&{7}\\{2}&{9}&{11, 13}\\{3}&{12}\\{6}\\{8,10}}.\]
Here ${\tt Heckeshape}(w)=(4,3,2,1,1)$ and indeed the length of the
first row of this shape equals ${\tt LIS}(w)=4$, whereas the 
length of the
first column equals ${\tt LDS}(w)=5$.
\end{Example}

\subsection{Proof of Proposition~\ref{prop:central}}
The claim that
$\mu_{n,q}$ is a probability distribution follows if {\tt Hecke} extends to 
provide a bijection between: 
\[W_{n,q} \mbox{ \ \ \ \ and \ \ \ \ } \Gamma_{n,q}:=\bigcup_{\lambda} INC(\lambda,q)\times
{\tt SsetVT}(\lambda,n),\]
where $\lambda$ satisfies (\ref{eqn:ph3}).

Associate to each word $w\in W_{n,q}$ the pair $(w,123\cdots n)$. Clearly
${\tt Hecke}$ injectively maps these pairs into $\Gamma_{n,q}$. 

To prove surjectivity, let $(P,Q)\in \Gamma_{n,q}$. Then
under ${\tt Hecke}^{-1}$, $(P,Q)$ corresponds to some pair
$(w,i)$. Now $i=123\cdots n$ since that is the only possible sequence
that can arise from a standard tableau $Q$. Also, since 
$W(w)=W({\tt word}(P))$, 
$w$ must use some subset
of $\{1,2,\ldots,q\}$. Thus $w\in W_{n,q}$. Hence 
$W_{n,q}\twoheadrightarrow \Gamma_{n,q}$. The claim (\ref{eqn:ph3})
is then clear from the properties of {\tt Hecke}.

Finally, from the above discussion it is immediate that {\tt Heckeshape}
is a sampling algorithm
for $\mu_{n,q}$. The bottleneck of the algorithm is the 
insertion process (a random uniform 
word $w\in W_{n,q}$ can be generated in $O(n\log(q))$ 
time). By (\ref{eqn:ph3}) we know that each
of the $n$ insertions demand at most ${q+1 \choose 2}$ operations.
Hence $O(nq^2)$ operations are needed. 
This completes the proof of Proposition~\ref{prop:central}.\qed

\section{Some further consequences of Theorem~\ref{thm:hard}}

\subsection{A generalization of the Erd\H{o}s-Szekeres theorem}

The following classic result is due to Erd\H{o}s-Szekeres \cite{Erdos}:

\begin{Theorem}
\label{thm:erdos}
Let $a,b\geq 1$. If $w\in S_{ab+1}$ then ${\tt LIS}(w)>a$  or
${\tt LDS}(w)>b$.
\end{Theorem}

It is known that this result can be readily deduced from Schensted's
results, see, e.g., \cite[Section~2]{Stanley:inc}.
Theorem~\ref{thm:hard} similarly leads to 
an extension of Theorem~\ref{thm:erdos} that relates {\tt LIS} and
{\tt LDS} to Coxeter length.

\begin{Proposition}
\label{prop:erdosgen}
Let $w\in W_{n,q}$. Suppose $1\leq a,b< q$ and 
\begin{equation}
\label{eqn:erdosgencond}
\ell(W(w))>\sum_{i=1}^{a}\min(b,q-i), \mbox{\ \ \ or equivalently\ \ \ }
\ell(W(w))>\sum_{j=1}^{b}\min(a,q-j),
\end{equation}
then ${\tt LIS}(w)>a$ or ${\tt LDS}(w)>b$; recall $W(w)$ is 
the permutation identified with $w$.
\end{Proposition}
\begin{proof}
If ${\tt LIS}(w)\leq a$ and ${\tt LDS}(w)\leq b$ then by 
Proposition~\ref{prop:central} and Theorem~\ref{thm:hard}, we have:
\[{\tt Heckeshape}(w)\subseteq (a\times b) \cap (q,q-1,\ldots,3,2,1).\]
Thus 
\[|{\tt Heckeshape}(w)|\leq \sum_{i=1}^{a}\min(b,q-i) = 
\sum_{j=1}^{b}\min(a,q-j).\]
Since $\ell(W(w))\leq |{\tt Heckeshape}(w)|$, the result then follows.
\end{proof}

\begin{Example}
If $q=4$ and $a=b=3$ then if $\ell(W(w))>8$ then ${\tt LIS}(w)>3$
or ${\tt LDS}(w)>3$. This inequality is tight in the sense that the
bound $8$ cannot be reduced: 
consider $w=2 \ 1 \ 3\ 4 \ 2 \ 3 \ 1 \ 2$. This is already a reduced
word of Coxeter length $8$, viewed as an element of $S_5$, and 
${\tt LIS}(w)={\tt LDS}(w)=3$. 
\end{Example}

Proposition~\ref{prop:erdosgen} generalizes Theorem~\ref{thm:erdos} because 
if $w\in S_{ab+1}$ is viewed as a Hecke word, we
have $\ell(w)=ab+1$ (any word where all letters are distinct is automatically reduced). Then set $q=ab+1$ and thus (\ref{eqn:erdosgencond}) is satisfied.

\subsection{Patience sorting for decks with repeated values}

In \cite{Aldous.Diaconis}, the Schensted correspondence was connected
to the one-person (solitaire) card game {\bf patience sorting}. 
We include a generalization of
this connection, which in
particular is a refinement of the {\tt LIS} claim of Theorem~\ref{thm:hard}.

In this game, a deck of cards labeled
$1,2,\ldots,n$ is shuffled and the cards are turned up one at a time
and dealt into piles on the table: a lower card
may be placed on top of a higher card, or put into  a new pile to the
right of the existing piles.
The goal of the game is to finish with as few piles as possible. 

For example, if $n=10$ and the deck is shuffled in the order
\[8 \ \  2 \ \   6 \ \ 3 \ \  4 \ \ 1\ \ 7 \ \ 10\ \ 9\]
then the top card $8$ is dealt onto the table. The $2$ can either be
placed to the right of the $8$ or on top of it -- suppose we chose the
latter scenario. Next the $6$ must be placed to the right of the pile
containing the $2$ and $8$, starting a new pile. At this stage, we have
$\begin{array}{cc}
2 & \\
8 & 6
\end{array}$.

The {\bf greedy strategy} is to always place the new card in the leftmost
pile possible. If we complete the game using this strategy, we would
obtain, successively:
\[\begin{array}{cc}
2 & 3\\ 
8 & 6
\end{array}
\mapsto 
\begin{array}{ccc}
2 & 3 & \\
8 & 6 & 4
\end{array}
\mapsto
\begin{array}{ccc}
1 &   & \\
2 & 3 & \\
8 & 6 & 4
\end{array}
\mapsto
\begin{array}{cccc}
1 &   &  & \\
2 & 3 & & \\
8 & 6 & 4 & 7
\end{array}
\mapsto 
\begin{array}{ccccc}
1 &   &  & & \\
2 & 3 & & & \\
8 & 6 & 4 & 7 & 10 \\
\end{array}
\mapsto 
\begin{array}{ccccc}
1 &   &  & & \\
2 & 3 & & & 9 \\
8 & 6 & 4 & 7 & 10 \\
\end{array}
\]
It is easy to prove that the top cards
increase from left to right throughout the game; Mallows \cite{Mallows:II}
and later independently Hammersley \cite[p.~362]{Hammersley:I} observed
that the number of piles at the end equals ${\tt LIS}(w)$, where
$w\in S_n$ is the permutation defining the shuffled deck. Finally, 
Aldous-Diaconis note that the first row of the insertion tableau under
Robinson-Schensted agrees with the top cards.

Aldous-Diaconis \cite[Section~2.4]{Aldous.Diaconis}) consider  
two variants of patience sorting where the deck has repeated entries, i.e.,
where all cards of the same rank 
(e.g., all Jacks) are equal. 
The two rules they consider are ``ties forbidden'' and ``ties allowed'', 
depending on whether or not a Jack can be placed on top of another Jack. 
They provide an analysis of the former case, relating it to the
Robinson-Schensted-Knuth correspondence. 

For example, if the shuffled deck is given by
$w= 2 \ \ 1 \ \ 4 \ \ 1 \ \ 3 \ \ 5 \ \ 3 \ \ 2 \ \ 5 \ \ 1 \ \  4 \ \ 
2\in W_{12,5}$,
then the result of playing patience (using the greedy strategy)
with ties forbidden
and allowed, respectively, are:
\[\begin{array}{ccccc}
  &   & 1 & 2 & \\
1 & 1 & 2 & 3 & 4\\
2 & 4 & 3 & 5 & 5
\end{array}
\mbox{ \ \ \ and \ \ \ }
\begin{array}{ccc}
  & 2 &  \\
1 & 2 &  \\
1 & 3 & 4  \\
1 & 3 & 5  \\
2 & 4 & 5
\end{array}.
\]

\begin{Proposition}
\label{prop:patience} Assume patience 
sorting is played with ties allowed, 
on a deck of $n$ cards with $q$ distinct types of cards (viewed as a word
$w\in W_{n,q}$). Then
\begin{itemize}
\item[(I)] The top cards of each pile at the termination of the game,
using the greedy strategy, as read from left to right, agree with the top 
row of the insertion tableau of ${\tt Hecke}(w)$.
\item[(II)] The optimal strategy 
(minimizing the number of piles created) is the greedy strategy, and
${\tt LIS}(w)$ piles are created.
\end{itemize}
\end{Proposition}
\begin{proof}
The proof of (I) is 
an easy induction, comparing the description of {\tt Hecke} with
the ``tied allowed'' rules of patience sorting.

For (II), by (I) and Theorem~\ref{thm:hard} 
we know the ``greedy strategy'' terminates with ${\tt LIS}(w)$ piles.
The proof of optimality is the same \emph{mutatis mutandis} as the one for
the original variant, see 
\cite[Lemma~1]{Aldous.Diaconis}.
\end{proof}

Briefly, probabilistic and statistical analysis on the 
``tied allowed'' case of patience sorting is possible, in analogy
to the work of \cite{Aldous.Diaconis}. 
Below we have tabulated the results of a Monte Carlo simulation
with $100,000$ trials on a standard $52$-card deck. 
Typically, the number of piles is 
between $8$ and $11$. The average number of piles is $9.2$
which, naturally, is less than the average number of piles
when the deck is totally ordered, which is $11.6$ as reported by
Aldous-Diaconis. So a deck ordering is
``lucky'' if the number of piles is less than $7$, which occurs
only about $3\%$ of the time.

\begin{table}[h]
\begin{tabular}{|l||l|l|l|l|l|l|l|l|}
\hline
number of piles &  6 & 7 & 8 & 9 & 10 & 11 & 12 & 13\\ \hline
frequency & 82 & 2993 & 20336 & 39039 & 27843 & 8489 & 1166 & 52\\ 
\hline
\end{tabular}
\caption{Monte Carlo simulation for standard 52-card deck
with $100,000$ trials. The average number of piles is
$9.2$.  }
\end{table}

\begin{figure}[h]
  \centering
  \epsfig{file=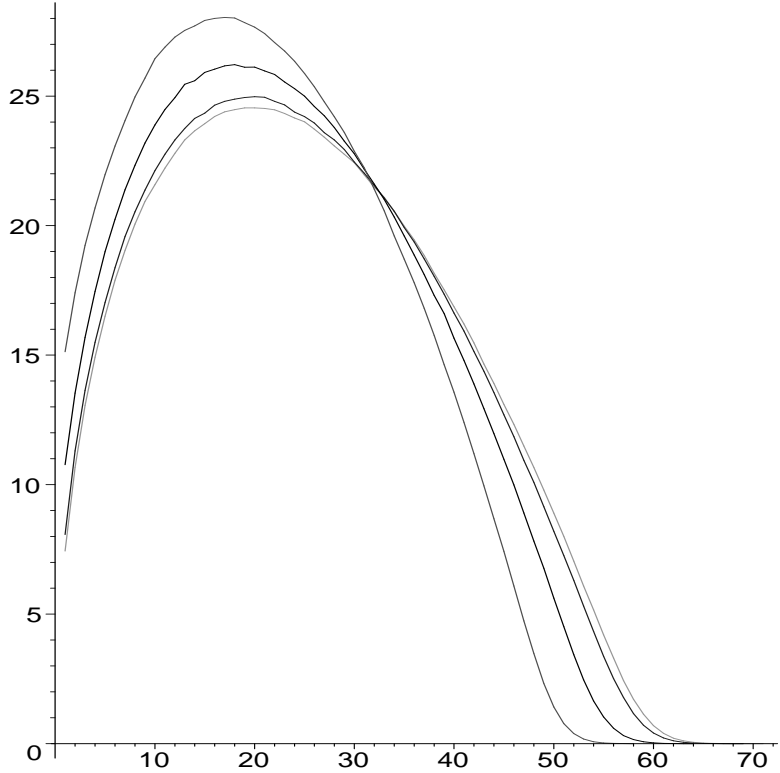,height=4in, width=4in}
\caption{\label{fig:pilesize} Simulation of mean pile sizes for $n=1000$ 
and $q=100,200, 600, 1000$ with
$10^4$ Monte Carlo trials. 
The $x$-axis indicates the position of the pile, counting
from the left.}
\end{figure}

Figure~\ref{fig:pilesize} shows that there
is a definite shape describing the mean pile sizes as $\alpha$ 
varies. Questions about the structure of this shape may be interpreted
as enriched questions related to {\tt LIS}. 
One can analyze such questions
using the dichotomy of Section~1.4; we do not pursue this here.

\section{Increasing tableau theory and the proof of Theorem~\ref{thm:hard}}
Now we show that the size of the first row of 
$\lambda={\tt Heckeshape}(w)$ 
computes ${\tt LIS}(w)$. Let $r(w,t)$ be the largest
index such that the longest strictly increasing subsequence 
ending at $w_{r(w,t)}$ has length $t$. For example, if
$w=2 \ 3\ 4 \ 1 \ 5 \ 2$ then $r(w,1)=4$, $r(w,2)=6$ and 
$r(w,3)=3$. We now in fact prove the following claim, which is
stronger than the {\tt LIS} assertion of Theorem~\ref{thm:hard}
(cf.~\cite[Prop.~7.23.10]{Stanley:text}):

\begin{Proposition}
\label{prop:refined_firstrow}
Suppose $w\in W_{n,q}$, $s={\tt LIS}(w)$ and $P$ is the 
insertion tableau of ${\tt Hecke}(w)$. Then the first row
of $P$ is given by $w_{r(w,1)},w_{r(w,2)},\ldots,w_{r(w,s)}$.
\end{Proposition}
\begin{proof}
By induction on $n$. The base case $n=1$ is trivial.

Suppose that the claim holds for $w^{\circ}:=w_1 w_2\cdots w_{n-1}$. Thus
if $P^{\circ}$ is the insertion tableau of ${\tt Hecke}(w^{\circ})$ then 
by induction, the first row is given by
\[w^{\circ}_{r(w^{\circ},1)}, w^{\circ}_{r(w^{\circ},2)},\ldots,
w^{\circ}_{r(w^{\circ},s^{\circ})},\]
where $s^{\circ}={\tt LIS}(w^{\circ})$. 

We consider the possibilities of what happens as we insert $w_n$ 
into the first row of $P^{\circ}$. We prove the desired 
conclusion holds for the case that 
$w^{\circ}_{r(w^{\circ},t)}< w_n$ for some maximally chosen
$t<s^{\circ}$; other cases
are similar. So after inserting $w_n$, the first row of $P$ is
\[w^{\circ}_{r(w^{\circ},1)}, w^{\circ}_{r(w^{\circ},2)},\ldots,
w^{\circ}_{r(w^{\circ},t)}, w_n, w^{\circ}_{r(w^{\circ},t+2)},\ldots,
w^{\circ}_{r(w^{\circ},s^{\circ})}.\]
The assumption shows that
the longest increasing subsequence $\omega$ 
in $w$ using $w_n$ is of length at least 
$t+1$, since we can adjoin $w_n$ to the length $t$ subsequence ending at
$w^{\circ}_{r(w^{\circ},t)}$. On the other hand, $\omega$ is of length
at most $t+1$ since if it were, say, of length $t+2$, there must be a
length $t+1$ increasing subsequence in $w^{\circ}$ ending at $w_r$
with $w_r<w_n$ and $r<r(w^{\circ},t+1)$. But this is a contradiction
of the definition of $r(w^{\circ},t+1)$ since we would then have a length
$t+2$ increasing subsequence ending at $w^{\circ}_{r(w^{\circ},t+1)}
(\geq w_n)$.

Since we have just shown $w_n=w_{r(w,t+1)}$, it is now clear that
$w_{r(w,h)}=w^{\circ}_{r(w^{\circ},h)}$ for $h\neq t+1$.  Thus,
the first row of $P$ satisfies the desired claim, and the induction
follows.
\end{proof}

In order to prove that the first column of $\lambda$ computes
${\tt LDS}(w)$, we need to draw a connection to 
\cite{Thomas.Yong:V}, where we developed a $K$-theoretic 
\emph{jeu de taquin} theory. Rather than repeat the setup in full here,
for brevity, we refer the reader to that paper
for the complete
background on $K$-{\tt rectification} and $K$-{\tt infusion} used below.

Although 
what follows also constitutes a proof of our {\tt LIS} claims,
we felt that including the direct proof via the stronger claim of 
Proposition~\ref{prop:refined_firstrow} was worthwhile.
However, a similarly direct proof of our {\tt LDS} claim seems harder.

Let 
\[\gamma_{n}=(n,n-1,n-2,\ldots,3,2,1)\] 
be the staircase shape. Also,
let $\lambda_{{\rm perm}(n)}=\gamma_{n}/\gamma_{n-1}$ be the {\bf permutation
shape} consisting of $n$ single boxes arranged along an antidiagonal.
Given $w\in W_{n,q}$, define $T_{w}\in {\tt INC}(\lambda_{{\rm perm}(n)})$
to be the tableau where $w_1,w_2,\ldots,w_n$ is arranged from
southwest to northeast. Also let $S\in {\tt SYT}(\gamma_{n-1})$ be the
{\bf superstandard} Young tableau, i.e., the one whose first row is labeled
$1,2,\ldots,n-1$, the second row is labeled by $n,n+1,n+2,\ldots,2n-3$ etc.
This latter tableau determines a \emph{particular} 
$K$-{\tt rectification} of $T_{w}$,
which by definition is $K$-${\tt infusion}_1(S, T_w)$, see 
\cite[Section~3]{Thomas.Yong:V}. (An important subtlety in $K$-theoretic
jeu de taquin is that $K$-{\tt rectification} 
depends on the order in which it
is performed, unlike the rectification of classical jeu de taquin. 
However, the order defined by $S$ is particularly nice.)

The following result is an analogue of a classical result linking the
Robinson-Schensted algorithm to the (ordinary) rectification of
$T_w$:

\begin{Theorem}
\label{thm:hecke_is_rect}
Let $w\in W_{n,q}$. Then
$K$-${\tt infusion}_1(S, T_w)$ is the insertion tableau of ${\tt Hecke}(w)$.
\end{Theorem}

\begin{proof}
We induct on $n$. The base cases $n=1,2$ are easy. We may assume that
the steps of the $K\mbox{-}{\tt infusion}_1$ that are defined by the ``inner'' labels
$\underline{n},\underline{n+1},\underline{n+2},\ldots, 
\underline{{n \choose 2}}$ results in a skew shape of the form
$P^{\circ}\star \tableau{w_n}$, as depicted below:
\begin{equation}
\label{eqn:depiction}
\ktableau{\underline 1 &\underline 2 & \cdots & {\underline{n-1}} &{w_n}\\{P^{\circ}_{1,1} }&{P^{\circ}_{1,2} }&{ \cdots }&{P_{1,n-1}^{\circ}}\\{P^{\circ}_{2,1} }&{\cdots }&{P_{2,n-2}^{\circ}}\\{\cdots }&{P_{n-2,2}^{\circ}}\\{P_{n-1,1}^{\circ}}}.
\end{equation}
The induction hypothesis is that $P^{\circ}$ is the insertion tableau
obtained by {\tt Hecke} inserting $w_1 w_2\cdots w_{n-1}$. 
(In the depiction of $P^{\circ}$ from (\ref{eqn:depiction}), 
some of the boxes with labels $P_{i,j}^{\circ}$ may be empty.) 
The non-underlined labels occupy the boxes of $P^{\circ}\star \tableau{w_n}$, 
whereas the underlined labels dictate the remaining steps
to perform to complete the $K$-${\tt infusion}_1$ computation (these steps
are recalled below).  

Hence it remains to show that the tableau obtained by the 
{\tt Hecke}-insertion $P^\circ\leftarrow w_n$ is the same as carrying out
the $K$-${\tt infusion}_1$ indicated by (\ref{eqn:depiction}), i.e.,
the operation
\begin{equation}
\label{eqn:Kinf_ind}
K\textrm{-}{\tt infusion}_1\left(\ktableau{{\underline 1} & {\underline 2}&\cdots&{\underline{n-1}}}, P^{\circ}\star\tableau{w_n}\right).
\end{equation} 
To do this, we first develop a technical fact.
In \cite[Section~1.1]{Thomas.Yong:V}, we defined the procedure {\tt switch},
which we restate now (in a more convenient form). Let 
${\tt Mixedtab}(\alpha,p,q)$ be the set of {\bf mixed tableaux}, which,
by definition, are tableaux of shape $\alpha$, 
each of whose boxes is filled with an entry
from one of two alphabets, 
$\{\underline
1,\dots,\underline p\}$ and $\{1,\dots,q\}$, such that, within each row
or column, the entries for each alphabet appear at most once. (No 
increasingness condition is demanded.) We also include the {\bf null tableau} 
$\emptyset$, as a special element of ${\tt Mixedtab}(\alpha,p,q)$.

Define an operator
\[{\tt switch}(\underline i,j):{\tt Mixedtab}(\alpha,p,q)\to
{\tt Mixedtab}(\alpha,p,q)\] 
as follows. Given $\emptyset\neq T\in {\tt Mixedtab}(\alpha,p,q)$, consider the subshape 
$S$ of $T$ consisting of boxes whose entry is either $\underline i$ or $j$.  
For each non-singleton connected component of $S$, interchange the 
$\underline i$'s and the $j$'s. If this results in a (non null) 
mixed tableau, then the result is that tableau. Otherwise the result is
$\emptyset$. By definition ${\tt switch}(\underline i,j)(\emptyset)=\emptyset$.

\begin{Example}
\label{eqn:a_switch}
Let $\alpha=(4,3,1)$ and $p=q=3$. Then $T\in {\tt Mixedtab}(\alpha,p,q)$
is given below, together with two different {\tt switch} computations applied
to it:
\[T\in\tableau{{\underline 2}&{\underline 1}&{\underline 3}&1\\
{\underline 1}&3&1\\2}, \  \ \ 
{\tt switch}({\underline 1}, 2)(T)
=\tableau{{\underline 2}&{\underline 1}&{\underline 3}&1\\2&3&1\\{\underline1}}, \ \ \
{\tt switch}({\underline 3},1)(T)=
\tableau{{\underline 2}&{\underline 1}&1&{\underline 3}\\ {\underline 1}&3
&{\underline 3}\\2}.
\]
On the other hand,
\[{\tt switch}(1, {\underline 2})\left(\tableau{1&{\underline 2}\\{\underline 1}\\{\underline 2}}\right)=\emptyset.\]
\end{Example}

The following lemma is easy to verify from the definitions:
\begin{Lemma} 
\label{lemma:viable}
If $i\ne j$ and $r\ne s$
then the operators ${\tt switch}(\underline i,r)$ and ${\tt switch}(\underline j,s)$ 
commute, i.e.,
\[{\tt switch}(\underline i,r) \ {\tt switch}(\underline j,s) \equiv
{\tt switch}(\underline j,s) \ {\tt switch}(\underline i,r)\]
is a relation in the algebra generated 
by {\tt switch} operators on ${\tt Mixedtab}(\alpha,p,q)$.
\end{Lemma}

The procedure described in \cite[Section~3]{Thomas.Yong:V} for
computing $K$-${\tt infusion}_1(A,B)$ is to consider the entries of $A$ as
being underlined, where the maximum entry of $A$ is $\underline p$, say, and 
the entries of $B$ as not underlined, where the maximum entry of $B$ is
 $q$.  Now perform the following sequence of ${\tt switch}$ 
operations, from left to right, as indexed by: 
\begin{equation}
\label{eqn:std}
(\underline p,1), (\underline p,2),\dots,
(\underline p,q), (\underline{p-1},1),\dots,
(\underline{p-1},q),\ldots,\ldots,(\underline 1,1),\ldots,
(\underline 1,q),
\end{equation}
We refer to this sequence of pairs (interchangeably, the corresponding
sequence of {\tt switch} operators) as the {\bf standard {\tt switch} 
sequence}.

The technical fact we need is that $K$-${\tt infusion}$ can
in fact be computed differently: 
a {\tt switch} sequence is called {\bf viable} if it is a ``shuffling'' of
(\ref{eqn:std}), in the following sense: 
\begin{itemize}
\item every $(\underline i,j)$ 
occurs exactly once, for $1\leq i\leq p$ and $1\leq j\leq q$;
\item  for any $1\leq i\leq p$, 
the pairs $(\underline i,1),\dots,(\underline i,q)$ occur
in that relative order; and 
\item for any $1\leq j\leq q$ the pairs $(\underline p,j),\dots
(\underline 1,j)$ occur in that relative order.  
\end{itemize}
This definition is explained by the proof of the following proposition:
\begin{Proposition}
\label{prop:viable}
Any viable {\tt switch} sequence can be used to calculate 
$K$-{\tt infusion}.
\end{Proposition}
\begin{proof} 
It is straightforward to show that one 
can obtain any viable {\tt switch} sequence from
the standard {\tt switch} sequence (\ref{eqn:std}) by repeated 
applications of the 
commutation relation of Lemma~\ref{lemma:viable}. 
\end{proof}

Thus, in view of Proposition~\ref{prop:viable}, 
to complete the induction
it suffices to construct a viable {\tt switch} sequence whose result
is the same as $P^{\circ}\leftarrow w_n$.
(A caution: in \cite{Thomas.Yong:V} it was shown that the standard {\tt switch}
sequence necessarily maintains increasingness along rows and columns of the
members of each alphabet. We will not prove that a
viable {\tt switch} sequence also achieves 
this during the intermediate steps of a $K$-{\tt infusion}. However, this 
does not play a logical role in how we apply 
Proposition~\ref{prop:viable} below.)

Let $y_1:=w_n$, and for $i>1$, let 
$y_i$ be the number which is inserted in row $i$ 
according to {\tt Hecke} insertion of $w_n$ into $P^{\circ}$.  
During a bumping step of {\tt Hecke} insertion, 
let $z_i$ be the smallest number already
in row $i$ which is greater that $y_i$. 

We say that a mixed tableau, obtained
after some number of {\tt switch} operations applied to 
$P^{\circ}\star\tableau{w_n}$, is in {\bf row $i$ normal form} if:
\begin{itemize}
\item the $i$-th row is of the form 
\[\ktableau{{\underline 1}&{\underline 2}&\cdots&{\underline{k-1}}&
y_i & {\underline{k}}&\cdots& {\underline t}},\] 
where $y_i$ has not yet moved from its
initial position in column $k$, or
\[\ktableau{{\underline 1}&{\underline 2}&\cdots&{\underline{k-1}}&
{\underline k} & {\underline{k+1}}&\cdots& {\underline t}},\] 
depending, respectively, on whether the {\tt Hecke} insertion $P^{\circ}\leftarrow w_n$ terminates
at row $i$ or after, or strictly earlier. Here the $i$-th row of $P^{\circ}$
has length $t$; and
\item all non-underlined symbols in rows $i+1$ and below have not
moved from their initial positions.
\end{itemize}

Having explained our general strategy, what remains is some tedious
but straightforward case
analysis to describe the viable {\tt switch} sequence we use: 
Our initial mixed tableau is in row~$1$ normal form. Now,
suppose we have arrived, after some sequence of applications
of (A), (B) and (C) below, at a mixed tableau in row $i$ normal form. 
There are three possibilities for a bumping step of
{\tt Hecke} insertion of $y_i$.

(A) \emph{$y_i$ is inserted, bumping $z_i$:}  
Consider the {\tt switch} sequence that moves $y_i$ 
to the left along row $i$: specifically, using 
\begin{equation}
\label{eqn:caseA_startmoves}
{\tt switch}({\underline{k-1}}, y_i), \ {\tt switch}({\underline{k-2}}, y_i),\
 {\tt switch}({\underline{k-3}}, y_i), \ldots,
\end{equation}
until it is directly above the $z_i$ in row $i+1$, 
and then, starting from the right, 
swaps each box in row $i$ having an underlined label with the one
directly below, which has a non-underlined label (this can always be done
since no label numerically equal to $y_i$ appears among the latter boxes, by assumption).  The
result is therefore in row $i+1$ normal form since $z_i$ doesn't move in this
process, it is the unique box with a non-underlined label in row $i+1$, and
$y_{i+1}=z_i$, as demanded.

Note that the non-underlined labels of the $i$-th row of the mixed tableau
we obtain after this process is the same as the $i$-th row of $P^{\circ}$,
with $z_i$ replaced by $y_i$, as desired. 

\begin{Example}
If $i=1$ and we started with 
\[\tableau{{\underline 1}&{\underline 2}&{\underline 3}&{\underline 4}&3\\
1&2&4&5\\2&3&5&6}\]
so that $y_i=3$ and $z_i=4$, then we begin by moving $y_i$ above $z_i=y_{i+1}$:
\[\tableau{{\underline 1}&{\underline 2}&3&{\underline 3}&{\underline 4}\\
1&2&4&5\\2&3&5&6}\]
We conclude with ${\tt switch}({\underline 3}, 5)$, 
${\tt switch}({\underline 2}, 2)$ and then finally ${\tt switch}({\underline 1},1)$, resulting in
\[\tableau{1&2&3&5&{\underline 4}\\{\underline 1}&{\underline 2}&4&{\underline 3}\\2&3&5&6},\]
which is in row $i+1$($=2$) normal form. Moreover, the non-underlined labels in 
row $i=1$ of this mixed tableau, namely $\tableau{1&2&3&5}$ agrees with the first row of  
$\tableau{1&2&4&5}\leftarrow y_i$, as desired.
\end{Example}

(B) \emph{$y_i$ is not inserted because of a horizontal violation:}
(That is, a label numerically equal to 
$y_i$ already appears in the row that we are {\tt Hecke} inserting
into.) Proceed as in (A) by moving 
$y_i$ to the left, until it is directly above $z_i$, i.e., apply 
(\ref{eqn:caseA_startmoves}). Now, 
``locally'' the situation in the column containing $z_i$ and the column
to its left is:
\[\cdots \tableau{\\  {\underline t} &y_i\\y_i&z_i} \cdots \mbox{ \ \ or \ \ }
\cdots\tableau{ &{\underline t}\\ {\underline t} &y_i\\y_i&z_i} \cdots, \]
for some $t$. The former case shows rows $i$ and $i+1$ when ${\underline t}$
does not appear in row $i-1$, 
while the latter
case also includes row $i-1$ if it does.
To the right of the column containing
$z_i$, we swap boxes with underlined and non-underlined labels,
as in (A). Then we perform the transformation
\begin{equation}
\label{eqn:case2_trans}
\cdots \tableau{ \\  {\underline t} &y_i\\y_i&z_i} \cdots \mapsto
\cdots \tableau{ \\ y_i &{\underline t}\\{\underline t}&z_i} \cdots
\mbox{\ \  or \ \ }
\cdots \tableau{ &{\underline t}\\{\underline t} &y_i\\y_i&z_i} \cdots \mapsto
\cdots \tableau{ &y_i\\ y_i &{\underline t}\\{\underline t}&z_i} \cdots
\end{equation}
respectively. After this transformation, 
we complete by swapping, right to left, the boxes to the left of the
$y_i$, as in (A). The result is in the demanded row $i+1$ normal form. 

Note that unlike (A), row $i$ still has a box $c$ 
with an underlined label in it.
However, when we work on the row $i+1$ normal form, the reader can check that
our descriptions will force that box $c$ to be filled by $z_i$:
e.g., if we apply case (A) next, this will occur when we execute the
switches (\ref{eqn:caseA_startmoves}),
if we apply case (B) next, this will occur when we execute either the
switches (\ref{eqn:caseA_startmoves}) or (\ref{eqn:case2_trans}), 
whereas if we apply (C) next, the
replacement will occur during the switches (\ref{eqn:case3_trans}) below.
Hence, in the end, we find that the 
non-underlined labels in row $i$ are the same as the 
$i$-th
row of the {\tt Hecke} insertion of $w_n$ into $P^{\circ}$, as desired. 

\begin{Example}
If $i=1$ and we started with
\[\tableau{{\underline 1}&{\underline 2}&{\underline 3}&{\underline 4}&{2}\\
{1}&{2}&{3}&{4}\\{4}}\]
then moving the ``$2$'' in row $1$ to the left, as in (A), gives
\[\tableau{{\underline 1}&{\underline 2}&{2}&{\underline 3}&{\underline 4}\\
{1}&{2}&{3}&{4}\\{4}}.\]
The remaining swaps give
\[\tableau{{\underline 1}&{\underline 2}&{2}&{\underline 3}&{\underline 4}\\
{1}&{2}&{3}&{4}\\{4}}\mapsto
\tableau{{\underline 1}&{\underline 2}&{2}&{4}&{\underline 4}\\
{1}&{2}&{3}&{\underline 3}\\{4}}
\mapsto
\tableau{{\underline 1}&{2}&{\underline 2}&{4}&{\underline 4}\\
{1}&{\underline 2}&{3}&{\underline 3}\\{4}}\mapsto
\tableau{{1}&{2}&{\underline 2}&{4}&{\underline 4}\\
{\underline 1}&{\underline 2}&{3}&{\underline 3}\\{4}},
\]
and the latter is in row $i+1=2$ normal form. 
The mild complication of this case, as suggested above, is that row $1$ is not
yet the same as the first row of the insertion tableau of {\tt Hecke}.

However, as we begin to work on this row $2$ normal form, we start
by using ${\tt switch}({\underline 2},3)$ (beginning to move the $3$ in row
$2$ to the left), giving:
\[\tableau{{1}&{2}&{3}&{4}&{\underline 4}\\
{\underline 1}&{3}&{\underline 2}&{\underline 3}\\{4}}.
\]
Hence row $1$ now does agree with {\tt Hecke} insertion.
\end{Example}
 
(C) \emph{$y_i$ will not be inserted because of a vertical violation (and a 
horizontal violation does not occur):} While
{\tt Hecke} inserting $y_i$ into row $i$ of $P^{\circ}$, 
$z_{i}$ is directly below a label (numerically equal to)
$y_i$ that is in row $i-1$. This prevents one from 
replacing $z_i$ by the $y_i$ we are inserting. 
(Note that this case can only occur if $i>1$.) Moreover, since there is
no horizontal violation, the number immediately to the left of 
$z_i$, say $b$, satisfies $b<y_i$. Notice, in order for a vertical violation
to occur, the previous step of {\tt Hecke} inserting $y_{i-1}$ 
must have been an instance of 
case (B) or (C) (since we are {\tt Hecke} 
inserting a number into row $i$ which
also appears in row $i-1$). Hence the
box directly above and to the left of the $y_i$
contains an underlined label $\underline t$.

Now we begin by switching 
all the boxes with underlined labels in row $i$ to the right of $y_i$, with the non-underlined
labels directly below them. The remaining switches, which involve 
the rows $i-1$ to $i+1$, in the column of $y_i$ and the 
column to its immediate left, are as given as follows:
\begin{equation}
\label{eqn:case3_trans}
\cdots \tableau{ & \underline t \\ \underline t & y_i \\ b & z_i }\cdots
\mapsto
\cdots \tableau{ & \underline t \\ b&y_i \\ {\underline t}&z_i}\cdots
\mapsto \cdots\tableau { & y_i\\ b&\underline t \\ \underline t & z_i}\cdots.
\end{equation}
As in (B), we finish by completing a sequence of swaps
involving the columns to the left of the $b$.
The result of this process is a tableau in row $i+1$ normal form.  

As at the conclusion of (B), row $i$ of the resulting
row $i+1$ normal form 
mixed tableau, still contains an underlined label. 
As in that case, this label will
be switched with $z_i$, as desired, during the forthcoming 
{\tt switches}.

\begin{Example} 
\label{exa:Ccase}
To give an example of (C), the previous insertion must have
been of type (B) or (C), so consider the following example:

$$\tableau {\underline 1 & \underline 2 & \underline 3 & 3 \\
            1            & 3            &  5 \\
            2            & 4            & 6 }.$$

After the first step, an insertion of type (B), we reach row 2 normal form:

$$\tableau { 1           & 3        &\underline 2 & \underline 3\\
             \underline 1 & \underline 2 & 5 \\
             2            & 4            &6}.$$

As compelled by the conditions of a viable switch sequence, we switch
$\underline 2$ with $4$ before we switch $\underline 2$ with 5, and as 
a result, as described above, we get to row 3 normal form:

$$\tableau{ 1         &3           &5 &\underline 3\\
            2         &4           & \underline 2\\
            \underline 1 & \underline 2 & 6\\}.$$

The final result is:

$$\tableau{ 1 & 3&5&\underline 3\\
            2&4&6\\
            6&\underline 1 &\underline 2}.$$

Again, this tableau agrees with $P^{\circ}\leftarrow 3$.
\end{Example}

Note that after each of (A), (B) and (C), when we reach row
$i+1$ normal form, 
$z_{i+1}$ is necessarily in a column weakly to the left of $y_{i+1}$,
because 
$P^{\circ}$ is an increasing tableau and
$y_{i+1}$ together with the rows strictly below row $i+1$ 
are, at this point, still unaltered.
This observation guarantees that the switches as described above can
actually be executed, 
i.e., the above descriptions
are well-defined. 

Using similar analysis one can give {\tt switch} sequences for 
the terminating steps of {\tt Hecke} insertion, such that one
maintains row $i$ normal form for all $i\leq \ell+1$, 
where $\ell$ is the number of
rows of $P^{\circ}$.
We leave the straightforward details to the reader.

These constructions then show, by induction on the number of rows $\ell$, 
that we have a sequence of {\tt switch} operations transforming
$P^{\circ}\star \tableau{w_n}$ into $P^{\circ}\leftarrow w_n$.

\emph{Conclusion of the proof of Theorem~\ref{thm:hecke_is_rect}:} by the
fact that $P^{\circ}$ is an increasing tableau, and by the definition 
of normal form, it is easy to see that 
the sequence of ${\tt switch}$ operations used forms a 
viable sequence, after suitable insertions of any trivial
${\tt switch}({\underline i}, r)$ operators (that is to say, switch
operations which do not have any effect on the tableau).  
We can therefore apply Proposition~\ref{prop:viable} 
as we claimed earlier.
\end{proof}

\begin{Example}
Continuing Example~\ref{exa:Ccase}, the switch sequence we obtain, 
by following the descriptions of the cases (B) and (C) that are needed is:
\[({\underline 3}, 2), ({\underline 2}, 3), ({\underline 2}, 4),
({\underline 2}, 5), (\underline 2, 6), ({\underline 1}, 6).\]
This is not quite a viable sequence: 
although our constructions guarantee that it satisfies
the second and third conditions to be a viable sequence, it fails the
first, since, e.g., 
$({\underline 3}, 1)$ doesn't appear in the sequence, since this {\tt switch}
is never needed. However, clearly we can simply insert this 
trivial {\tt switch}, along with the others that are missing, giving
the viable sequence:
\[(\underline 3, 1), {\bf({\underline 3}, 2)}, (\underline 3,3),
(\underline 3, 4), (\underline 3,5), (\underline 3,6),
(\underline 2,1), (\underline 2,2),
{\bf ({\underline 2}, 3)}, 
{\bf ({\underline 2}, 4)},
{\bf ({\underline 2}, 5)}, {\bf (\underline 2, 6)},\] 
\[(\underline 1,1), (\underline 1,2), (\underline 1,3), (\underline 1,4),
(\underline 1,5),
{\bf ({\underline 1}, 6)}.
\]
The action of this viable sequence on the original mixed tableau is therefore
the same as the original {\tt switch} sequence, which we highlight 
in boldface. This viable sequence also happens to be the standard {\tt switch} 
sequence, although it needn't be in general. Hence 
\[K\mbox{-}{\tt infusion}_1\left(\tableau{1&2&3}, 
\tableau {& & & 3 \\
            1            & 3            &  5 \\
            2            & 4            & 6 }\right)=\tableau{1&3&5\\2&4&6}\leftarrow 3,\]
in agreement with Theorem~\ref{thm:hecke_is_rect}.
\end{Example}

\medskip
In \cite[Theorem~6.1]{Thomas.Yong:V} we showed that the first row of
$K$-${\tt infusion}_1(R, T_w)$ has length ${\tt LIS}(w)$, for any
increasing tableau $R$ of shape $\gamma_{n-1}$. So
by Theorem~\ref{thm:hecke_is_rect}, 
the first row of ${\tt Heckeshape}(w)=K\mbox{-}{\tt
infusion}_1(S, T_w)$ has length ${\tt LIS}(w)$. By symmetry,
\cite[Theorem~6.1]{Thomas.Yong:V} also implies that
the first column of $K$-${\tt infusion}_1(R, T_w)$ has length ${\tt
LDS}(w)$. Hence the {\tt LDS} claim follows. This completes the
proof of Theorem~\ref{thm:hard}.
\qed

Given $w=w_1 w_2 \cdots w_n$, define ${\tt rev}(w)=w_n w_{n-1} \cdots w_1$.
The following is symmetry statement is
immediate from Theorem~\ref{thm:main}, since 
${\tt LIS}(w)={\tt LDS}({\tt rev}(w))$:
\begin{Corollary}
\label{thm:partialsymm}
Let $\lambda={\tt Heckeshape}(w)$ and $\mu={\tt Heckeshape}({\tt rev}(w))$.
Then $\lambda_1=\mu_1'$ and $\mu_1=\lambda_1'$ where $\lambda'$ and $\mu'$
are the conjugate shapes of $\lambda$ and $\mu$ respectively.
\end{Corollary}

A warning is needed:
unlike Robinson-Schensted correspondence setting, with {\tt Hecke},
one \emph{cannot} conclude that the
insertion tableaux associated to $w=w_1 w_2 \cdots w_n$ 
and ${\tt rev}(w)=w_n w_{n-1} \cdots w_1$
differ only by a reflection across the main diagonal.
A counterexample is $w=1\ 3\ 4\ 2\ 2$. 
(In \cite{Schensted}, the symmetry property of the Robinson-Schensted 
correspondence, was applied to prove the {\tt LDS} claim in the 
classical version of Theorem~\ref{thm:hard}.)

\begin{Problem}
Give an explicit description of ${\tt Hecke}({\tt rev}(w))$
in terms of ${\tt Hecke}(w)$.
\end{Problem}

Finally, Greene \cite{Greene} has given an explanation of the other rows of the
shape $\lambda$ associated to a permutation $w\in S_n$ under the Robinson-Schensted 
correspondence: $\lambda_1+\ldots+\lambda_i$ equals the maximal size of
a union of $i$ disjoint increasing subsequences of $w$. 

However, we could not find any extension of 
Greene's theorem in the {\tt Hecke} context. The naive tries do not work:
Since $|\lambda|\leq n$, the simplest case
to analyze is when $|\lambda|=n$. The example $w=2 \ 1 \ 2 \ 3 \ 2$ corresponds
to $\lambda=(3,2)$; this shows that it is not valid to merely
replace ``increasing'' by ``strictly increasing'' in Greene's theorem, since
that would predict $\lambda=(3,1,1)$.

\section{Probabilistic combinatorics and proofs of Theorems~\ref{thm:desc} and~\ref{thm:main}} 

\noindent
\emph{Proof of Theorem~\ref{thm:desc}:}
 Let $w_0=\left(\begin{array}{ccccccc}
1&2&3&\ldots &q-1&q&q+1\\
q+1 & q & q-1 & \ldots & 3 & 2 & 1\end{array}\right)$ be 
the word in $S_{q+1}$ of maximal Coxeter length. Hence
$\ell(w_0)={q+1\choose 2}$. 
This is the unique permutation in $S_{q+1}$ with this length. 

We need the following
lemma, which characterizes when ${\tt Heckeshape}(w)$ is
maximized.
\begin{Lemma} 
\label{lemma:w0}
For $w\in W_{n,q}$, 
${\tt Heckeshape}(w)=(q,q-1,\ldots,3,2,1)$ if and only if
$W(w)=w_0$.
\end{Lemma}
\begin{proof}
First suppose $W(w)=w_0$.
Under ${\tt Hecke}(w)$, the insertion tableau $P$ satisfies \linebreak
$W({\tt word}(P))=w_0$. 
Hence the
shape of $P$ has at least ${q+1 \choose 2}$ boxes and 
so by Proposition~\ref{prop:central} it must be $(q,q-1,\ldots,3,2,1)$.

Conversely, if ${\tt Heckeshape}(w)=(q,q-1,\ldots,3,2,1)$, then
note that there is a unique increasing filling $P$ of that shape
(using $1,2\ldots, q$ in the first row, $2,3,\ldots,q$ in the
second row, etc). Then it is well-known that $W({\tt word}(P))=w_0$.
\end{proof}

In view of the Lemma~\ref{lemma:w0}, 
the Theorem will follow if we can show that
\begin{equation}
\label{eqn:almostsurely}
W(w)=w_0 \mbox{ \ \  almost surely, as $n\to\infty$}.
\end{equation}
(We conjecture this to be true whenever 
$0\leq \alpha< \alpha_{\rm critical}=\frac{1}{2}$. This would
imply Conjecture~\ref{conj:thebigone} for this entire range.)

Set $w(k):=w_1\cdots w_k$. Then either
\[\ell(W(w(k) \ w_{k+1}))=\ell(W(w(k)))+1 \mbox{\ \  or\ \  }
\ell(W(w(k)))\]
depending on whether the simple reflection $W(w_{k+1})$ 
(say equal to $s_t=(t \ t+1)$) is an {\bf ascent} of $W(w(w))$ or not.
(An {\bf ascent} occurs at a position $t$ for a permutation $\pi$ if  
$\pi(t)<\pi(t+1)$.)

Provided that $\pi\neq w_0$, $\pi$ has at least one ascent. Thus when
$W(w_k)\neq w_0$, the
probability that $w_{k+1}$'s introduction increases the Coxeter length is
at least $\frac{1}{q}$. 

Let 
\[E_k=\mbox{the event that $\ell(W(w(k)))<{q+1\choose 2}$}.\]
Related to this, let $Y_i\in\{0,1\}$ be Bernoulli distributed 
with parameter $\frac{1}{q}$. Set
\[Z_k:=Y_1+\cdots +Y_k.\]
Clearly,
\begin{equation}
\label{eqn:comparisonineq}
Prob\left(E_k\right)\leq 
Prob\left(Z_k<{q+1\choose 2}\right).
\end{equation}
We now show that when 
\[k=O(q^{3+\epsilon}), \mbox{\ \ for $0<\epsilon$},\]
the righthand side of the inequality (\ref{eqn:comparisonineq})
goes to zero as $q\to\infty$.

This is a simple application of (a special case of) Bennet's 
large deviation inequality, see, e.g., \cite[Cor.~2.4.7]{Dembo.Zeitouni}: 
suppose $X_i$ are independent, mean zero random variables with 
$|X_i|\leq 1$. Set
$S_{k}=\sum_{i=1}^k X_i$. Then for $y\geq 0$ we have
\begin{equation}
\label{eqn:bennet}
Prob(k^{-\frac{1}{2}}S_k\geq y)\leq e^{-\frac{y^2}{2}}.
\end{equation}
To apply this to our setting, let $X_i=-Y_i+\frac{1}{q}$. Hence
$S_k=-Z_k+\frac{k}{q}$. Then with $a=o(q^2)$
\begin{eqnarray}\nonumber
Prob(Z_k<a) & = & Prob\left(-S_k<a-\frac{k}{q}\right) = 
Prob\left(S_k\geq \frac{k}{q}-a\right)\\ \nonumber
& \leq & Prob\left(S_k> \frac{k}{2q}\right) \mbox{\ \ (for $q$ large, since $a=o(q^2)$)}\\ \nonumber
& = & Prob\left(k^{-\frac{1}{2}}S_k\geq \frac{\sqrt{k}}{2q}\right)
\leq e^{-\frac{k}{8q^2}}\to 0, \mbox{as $q\to\infty$}.\nonumber
\end{eqnarray}
The result then follows.
\qed

In the above argument, we interpreted $w\in W_{n,q}$ as
a random walk in $S_{q+1}$ that begins at 
the identity and works its way up in the weak Bruhat order
to $w_0$. At each step the probability of going up is 
\emph{at least} $\frac{1}{q}$ (as we have used), 
but is larger in general. However, since this probability varies,
even for permutations with the same Coxeter length, a more
refined analysis is needed to push the argument we have
used further, up towards $\alpha_{\rm critical}=\frac{1}{2}$. 

\noindent
\emph{Proof of Theorem~\ref{thm:main}:} 
For $0\leq \alpha<\alpha_{\rm critical}=\frac{1}{2}$ we will apply 
an argument similar to that for Theorem~\ref{thm:desc}.

Given $u\in W_{n,q}$ let 
\[m(u)=\max_{t\geq 1}\ 1,2,\ldots, t \mbox{\ \  is a subsequence
of $u$}.\]
Let $w(k)=w_1\cdots w_k$ and set
\[E_k =\mbox{the event that $m(w(k))<q$}.\]
Provided $E_k$ occurs, then
\[m(w(k)w_{k+1})=m(w(k))+1\]
with probability $\frac{1}{q}$, and is equal to 
$m(w(k))$ otherwise.

Let $\{Y_i\}$ and $Z_k=Y_1+\cdots +Y_k$ be discrete random variables, where
$Y_i$ is Bernoulli distributed with parameter $\frac{1}{q}$. Now,
\[Prob(E_k)=Prob(Z_k<q).\]
Thus it will be enough to show that when 
$k=O(q^{2+\epsilon}) \mbox{\  for $\epsilon>0$}$ then 
\[Prob(Z_k<O(q^{1+\epsilon}))\to 0 \mbox{\ as $q\to \infty$}.\] 
This is another application
of the large deviation inequality (\ref{eqn:bennet}). 
 
For $\alpha_{\rm critical}=\frac{1}{2}<\alpha\leq 1$ we use a proof provided for us by
O.~Zeitouni: $E({\tt LwIS})$, the expected
length of the longest weakly increasing subsequence of $w\in W_{n,q}$ 
(with $\alpha>\alpha_{\rm critical}=\frac{1}{2}$) is known to satisfy
\[E({\tt LwIS})\approx 2\sqrt n;\]
see \cite[Theorem~1.7]{Johansson}.
The argument shows that the difference between the {\tt LIS}
and {\tt LwIS} of $w$ is typically small.

Let ${\tt LwIS}_{a,b}$ be the random variable for the value of 
${\tt LwIS}(w)$ of a random uniform word
$w\in W_{\lfloor a\rfloor, \lfloor b\rfloor}$ where
$\lfloor a\rfloor$ is the integer part of $a$, etc. Similarly
define ${\tt LIS}_{a,b}$ where ${\tt LIS}$ replaces ${\tt LwIS}$.

Fix $\epsilon>0$ and let $L_0=L_0(\epsilon)$ be large enough
such that
\begin{equation}
\label{eqn:inflim}
{\rm inf}_{L>L_0}{\underline \lim}_{n\to\infty}Prob\left({\tt LwIS}_{L^2(1-\epsilon),\frac{qL}{\sqrt n}}>2(1-4\epsilon)L\right)>1-\epsilon.
\end{equation}
We need a ``graphical'' representation of a word in $W_{n,q}$: 
consider a $q\times n$
rectangle subdivided into unit squares. In each of the $n$ columns, one
places a single ``dot'' in one of the $q$ rows. The set of such 
configurations is in obvious bijection with words in $W_{n,q}$.

Given $L$, draw $\sqrt n/L$ smaller rectangles of dimension 
$\frac{qL}{\sqrt n}\times L{\sqrt n}$ along an antidiagonal
inside the $q\times n$ rectangle, as
depicted in Figure~\ref{fig:Zeit} below.

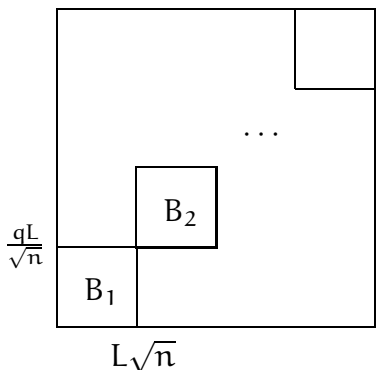
\begin{figure}[h]
\[\begin{picture}(120,120)
\put(0,0){\framebox(120,120)}
\put(30,0){\line(0,1){30}}
\put(10,10){$B_1$}
\put(20,-15){$L\sqrt n$}
\put(-20, 28){$\frac{qL}{\sqrt n}$}
\put(0,30){\line(1,0){30}}
\put(30,30){\framebox(30,30)}
\put(40,40){$B_2$}
\put(90,90){\line(1,0){30}}
\put(90,90){\line(0,1){30}}
\put(70,70){$\cdots$}
\end{picture}
\]
\caption{\label{fig:Zeit} $\alpha>\alpha_{\rm critical}=\frac{1}{2}$ case of proof of Theorem~\ref{thm:main}}
\end{figure}

Label the $i$-th southwest most box $B_i$. Let $N_i$ be the random variable
giving the number of dots in $B_i$. Notice that the $N_i$'s are independent.
Also, the dots inside $B_i$ define a word, and we can speak of 
${\tt LIS}(B_i)$ and ${\tt LwIS}(B_i)$, the length of the longest strictly
(respectively, weakly) increasing subsequence of that word.

Say that $B_i$ is {\bf good} if the following conditions simultaneously hold:
\begin{itemize}
\item[(a)] $N_i\geq L^2(1-\epsilon)$;
\item[(b)] ${\tt LwIS}(B_i)\geq 2(1-4\epsilon)L$; and
\item[(c)] no two dots in $B_i$ have the same height (hence ${\tt LwIS}(B_i)={\tt LIS}(B_i)$).
\end{itemize}

Now, we have 
\[E(N_i)=L{\sqrt n}\times\left(\frac{qL}{\sqrt n} \times \frac{1}{q}\right)=L^2\]
and we claim for an $L_1$ sufficiently large, for $L\geq L_1$, 
and for all $n$ large, we have:
\[Prob(N_i\leq L^2(1-\epsilon))\leq \epsilon.\]
The proof is a standard argument: let $Y_k$ be
the indicator random variable which evaluates to $1$ if column $k$ has
a dot that lies in the box $B_i$ that occupies that column (hence with
probability $L/\sqrt{n}$), and evaluates to $0$ otherwise. Hence
\begin{eqnarray}\nonumber
Prob(N_i\leq L^2(1-\epsilon)) & = & Prob(\sum_{k=1}^{L{\sqrt n}} Y_k
\leq L^2(1-\epsilon))\\ \nonumber
& = & Prob(\sum_{k=1}^{L{\sqrt n}}(Y_k-E(Y_k))<L^2(1-\epsilon)-L^2)\\ \nonumber
& = & Prob(\sum_{k=1}^{L{\sqrt n}}(Y_k-E(Y_k))<-L^2\epsilon)\\ \nonumber
& \leq & \frac{{L\sqrt n}\cdot (L/\sqrt{n})}{L^4\epsilon^2} = \frac{1}{L^2\epsilon^2},\nonumber
\end{eqnarray}
where the previous line is an application of Chebyshev's inequality.
Now take $L$ sufficiently large (bigger than some $L_1$)
so $\frac{1}{L^2 \epsilon^2}\leq \epsilon$.

Assuming $L\geq L_0$, if the event (a) occurs, then with probability
at least $1-\epsilon$, when $n$ is large, (b) holds, because of 
the definition of $L_0$ using (\ref{eqn:inflim}).

The probability of the event (c) not occurring is bounded above (using a
union bound) by
\begin{eqnarray}\nonumber
\#\mbox{rows}\times Prob(\mbox{two dots share a given height}) & \leq &
\frac{qL}{\sqrt n} \times (L\sqrt n)^2\times (1/q)^2\\ \nonumber
& = &
L^3\frac{\sqrt n}{q}\to 0,\nonumber
\end{eqnarray}
because $n=o(q^2)$.

So, 
\[Prob(B_{i}\mbox{ \ \ is good})\geq 1-3\epsilon
\]
for $L>\max(L_0, L_1)$ and $n$ large. 

Another standard argument with Chebyshev's inequality shows 
with high probability, say at least $1-\epsilon$, and for $n$ large, 
the number of good boxes is 
\[\frac{\sqrt n}{L}(1-\epsilon)(1-3\epsilon)\geq\frac{\sqrt n}{L}
(1-4\epsilon).\] 
Hence, with that probability, for $w\in W_{n,q}$
\[{\tt LIS}_{n,q}\geq \frac{\sqrt n}{L}(1-4\epsilon)\times 2(1-4\epsilon)L\geq
2{\sqrt n}(1-8\epsilon).\]
Since 
\[2{\sqrt n}(1-8\epsilon)\leq E({\tt LIS}_{n,q})\leq 
E({\tt LwIS}_{n,q})\approx 2\sqrt{n}\]  
the $\alpha>\alpha_{\rm critical}=\frac{1}{2}$ case follows by taking $\epsilon\to 0$, 
completing the proof of the theorem.
\qed

\section{Appendix (by A.~Yong and O.~Zeitouni)}

The goal of this appendix is to present a proof of the following result:

\begin{Theorem}
\label{thm:appendix}
Let $q=k\sqrt n+\mbox{lower order terms}$. Then
\[E({\tt LIS})\approx\beta(k){2\sqrt n}\] 
where
\begin{equation}\nonumber
\beta(k)=\left\{
\begin{array}{cc}
\frac{k}{2} & \mbox{if $0<k\leq 1$}\\
\frac{2-k^{-1}}{2} & \mbox{if $k>1$.}
\end{array}\right.
\end{equation}
\end{Theorem}
However, in order to prove this statement, we need to work with another
variant of Plancherel measure, utilized, e.g., by \cite{Biane} and alluded to 
in Section~1.5 of the main text. Our approach parallels the one developed
in \cite{Logan.Shepp, Vershik.Kerov} to prove $E({\tt LIS})=2\sqrt n$
in the permutation case, by utilizing work of \cite{Biane}. 

\subsection{Preliminaries}
A {\bf semistandard Young tableau} of shape $\lambda\in {\mathbb Y}$ with
labels from $\{1,2,\ldots, q\}$ is a filling of the Young shape $\lambda$
with these labels so that the entries weakly increase along rows, and 
strictly increase along columns. For example, if $\lambda=(2,1)$ and $q=2$
there are two such tableaux: $\tableau{1&1\\2}$ and $\tableau{1&2\\2}$.
Let $g^{\lambda}(q)$ denote the number of such tableaux.

Define the {\bf Plancherel-{\tt RSK} measure} $\nu_{n,q}$ on the set ${\mathbb Y}_n$
of Young diagrams $\lambda$ with $n$ boxes, by declaring that a random
Young shape $\lambda_{n,q}$ occurs with probability
\[Prob(\lambda_{n,q}=\lambda)=\frac{1}{q^n}f^{\lambda}g^{\lambda}(q).\]
We make no claims of originality in this definition. Indeed, this is the same
measure studied in, e.g., \cite{Biane}; 
although there the measure is defined
in terms of dimensions of irreducible $S_n$ and $GL_n({\mathbb C})$ modules
associated to $\lambda$; the equivalence is well-known. The fact that 
$\nu_{n,q}$ is in fact a probability distribution follows from either
Schur-Weyl duality, as in Section~1.5, or by the {\tt RSK} algorithm, see, 
e.g.,\cite[Section~7.11]{Stanley:text}.

A crucial advantage of $\nu_{n,q}$ for the purposes of
understanding $E({\tt LIS})$, in comparison to Plancherel-Hecke
measure, is that both $f^{\lambda}$ and $g^{\lambda}(q)$ have simple
multiplicative formulas. This makes it more readily analyzed using 
ideas of \cite{Logan.Shepp, Vershik.Kerov}, which we modify to the
present setting.

Given a box $u\in\lambda$, define the {\bf hook-length} associated to 
$u$ to be $H(u):=A(u)+L(u)+1$ where $A(u)$ is the number of boxes strictly to the right
of $u$, and in the same row, and $L(u)$ is the number of cells strictly
below $u$ and in the same column. Then we have 
\cite[Chapter~7]{Stanley:text}, the {\bf hook-length formula} and 
{\bf hook-content formulas}, respectively:
\[f^{\lambda}=\frac{n!}{\prod_{u\in\lambda}H(u)} \mbox{\ \ and \ \ }
g^{\lambda}(q)=\prod_{u\in\lambda}\frac{q+C(u)}{H(u)},\]
where in the second formula $C(u)$ is the {\bf content} of $u$, the
column index of $u$ minus the row index of $u$. So for example, if
$\lambda=(4,3,2)$, the contents are given by
\[\tableau{0&1&2&3\\-1&0&1\\-2&-1}.\]

\subsection{Plancherel-{\tt RSK} as a Markov measure} 
{\bf Young's lattice} is 
the poset structure on ${\mathbb Y}$ where $\lambda\leq \mu$
if the shape of $\lambda$ is contained in the shape of $\mu$. We write
$\lambda\to \mu$ to denote a covering relation in this poset, i.e., 
where $\mu$ is obtained from $\lambda$ by adding a single box at a corner.

Define a Markov process on ${\mathbb Y}$ with the transition probabilities
\[Prob(\lambda\to \mu)=\frac{g_{\mu}(q)}{qg_{\lambda}(q)}.\]
We need the following lemma, that in particular shows that 
Plancherel-{\tt RSK} measure is a Markov measure with the above 
transition probabilities. 

\begin{Lemma}
\begin{itemize}
\item[(I)] $\sum_{\Lambda:\lambda\to \mu} Prob(\lambda\to \mu)=1$
\item[(II)] $\sum_{\lambda:\lambda\to \mu} Prob(\lambda\to \mu)
\nu_{n,q}(\lambda)=\nu_{n+1,q}(\mu)$
\end{itemize}
\end{Lemma}
\begin{proof}
The claim (I) is equivalent to 
\[\sum_{\Lambda:\lambda\to \mu} g^{\mu}(q)=qg^{\lambda}(q).\]
This follows from the following \emph{Pieri rule} for Schur polynomials
\[\sum_{\Lambda:\lambda\to \mu} s_{\mu}(x_1,\ldots, x_q)=
s_{(1)}(x_1,\ldots,x_q)\cdot s_{\lambda}(x_1,\ldots,x_q).\]
See \cite[Theorem~7.15.7]{Stanley:text}. Here 
\[s_{\lambda}(x_1,\ldots,x_q)=\sum_{T}{\bf x}^T\]
is the {\bf Schur polynomial}, where the sum is over all semistandard
Young tableaux of shape $\lambda$ with entries from $\{1,2,\ldots,q\}$,
${\bf x}^T=x_1^{i_1}x_2^{i_2}\cdots x_q^{i_q}$, and 
$i_j$ is the number $j$'s used in $T$. In particular, (I) is immediate from
$g^{\lambda}(q)=s_{\lambda}(1,1,\ldots,1)$.

For (II), the claim is
\[\sum_{\lambda:\lambda\to \mu}
\left(\frac{g^{\mu}}{qg^{\lambda}}\right)\left(\frac{f^{\lambda}g^{\lambda}}{q^n}\right)=\frac{f^{\mu}g^{\mu}}{q^{n+1}},\]
that is,
$\sum_{\lambda:\lambda\to \mu} f^{\lambda}=f^{\mu}$, which is
well-known (and straightforward from the definitions).
\end{proof}

\subsection{Conclusion of Proof of Theorem~\ref{thm:appendix}}
Work of Biane \cite[Theorem~3]{Biane} describes the typical shape under 
Plancherel-{\tt RSK} after the rescaling 
$f\mapsto \frac{1}{2\sqrt n}f(2\sqrt n)$. Biane's theorem implies that
\begin{equation}
\label{eqn:leq}
E({\tt LIS})\geq \beta(k)
\end{equation} 
but not 
\begin{equation}
\label{eqn:geq}
E({\tt LIS})\leq \beta(k).
\end{equation}
Briefly, we explicate 
how his work applies to our situation (the 
reader is directed to the original source for details): 
Biane works with the coordinate 
axes rotated $45$-degrees 
counterclockwise, as in, e.g., \cite{Vershik.Kerov, Vershik.Kerov85}. 
His aforementioned theorem
states that if $\{f_n\}_{n=1}^{\infty}$ is any sequence of (rescaled and
rotated) Young diagrams, then for any $\epsilon>0$
we have
\begin{equation}
\label{eqn:Bianesresult}
\lim_{n\to\infty} Prob\left(\sup_{u\in {\mathbb R}}|f_{n}(u)-P_{\frac{1}{k}}(u)|>\epsilon\right)\to 0,
\end{equation}
where the probability is computed with respect to $\nu_{n,q}$. Also, 
$P_{\frac{1}{k}}$ is Biane's limit shape, which has the property that
it meets the line $y=x$ at a distance $\beta(k)$ from the origin. In other
words, the ``first column'' $C$ of $P_{\frac{1}{k}}$ satisfies
\begin{equation}
\label{eqn:C}
C=\beta(k).
\end{equation}
For each $n$, let $C_n$ be the length of the first column of $f_n$.
From (\ref{eqn:Bianesresult}) it follows that for any $\epsilon>0$,
\begin{equation}
\label{eqn:Cn}
\lim_{n\to\infty} Prob(C_n<C-\epsilon)\to 0.
\end{equation}
Moreover, since it is known that for any $w\in W_{n,q}$ the first column of the
Young diagram associated to ${\tt RSK}(w)$ equals ${\tt LDS}(w)$, 
(\ref{eqn:leq}) follows immediately from (\ref{eqn:C}) and (\ref{eqn:Cn})
combined. 

Note that the above argument does not also prove (\ref{eqn:geq}) since
(\ref{eqn:Bianesresult}) does not rule out the possibility that 
$\{f_n\}_{n=1}^{\infty}$ consists of Young diagrams with ``tails'' along
the $y=x$ axis that both ``lengthen'' and ``thin out'' as $n\to\infty$.
Therefore, it remains to verify (\ref{eqn:geq}).

To do this, we modify an argument found in \cite{Vershik.Kerov85}, 
which
establishes the analogous assertion in the permutation case: consider
the set ${\mathbb Y}^{\infty}$ of all sequences of Young diagrams 
\[{\bf \lambda}=(\lambda^{(1)},\lambda^{(2)}, \lambda^{(3)},\ldots,
\lambda^{(i)},\ldots)\]
where $\lambda^{(i)}\to \lambda^{(i+1)}$ for $i\geq 1$. 

For a Young diagram
$\lambda$, let $\lambda^{\downarrow}$ denote the diagram obtained by adding
a single box to $\lambda$, in the first column. For each integer $i\geq 1$,
define the indicator function $\psi_i:{\mathbb Y}^{\infty}\to \{0,1\}$ by setting
$\psi_i({\bf \lambda})=1$ if $\lambda^{(i)}=(\lambda^{(i-1)})^{\downarrow}$, 
and setting $\psi_i({\bf \lambda})=0$ otherwise. 

Studying the expectation of $\psi_i$ we have:
\begin{eqnarray}\nonumber
E(\psi_i)^2 & = & \left(\sum_{\lambda}\nu_{i-1,q}(\lambda)\cdot Prob(\lambda\to \lambda^{\downarrow})\right)^2\\ \nonumber
& \leq &  \sum_{\lambda}\nu_{i-1,q}(\lambda)\cdot Prob(\lambda\to \lambda^{\downarrow})^2 \mbox{\ \ \ \ \ \ \ \ (Cauchy-Schwarz inequality)}\\ \nonumber
& = & \sum_{\lambda} \nu_{i-1,q}(\lambda)
\frac{f^{\lambda^{\downarrow}}g^{\lambda^{\downarrow}}(q)}{f^{\lambda}g^{\lambda}(q)}\cdot \frac{1}{q} \cdot \frac{1}{q}\cdot \frac{g^{\lambda^{\downarrow}}(q)/g^{\lambda}(q)}{f^{\lambda^{\downarrow}}/f^{\lambda}}\\ \nonumber
& = & \frac{1}{q}\sum_{\lambda} \nu_{i,q}(\lambda^{\downarrow})\frac{g^{\lambda^{\downarrow}}(q)/f^{\lambda^{\downarrow}}}{g^{\lambda}(q)/f^{\lambda}},\nonumber
\end{eqnarray}
where we have just used 
\[\nu_{i,q}(\lambda^{\downarrow})=\nu_{i-1,q}(\lambda)
\frac{f^{\lambda^{\downarrow}}g^{\lambda^{\downarrow}}(q)}{f^{\lambda}g^{\lambda}(q)}\cdot \frac{1}{q}.\]

Let
\[L(\lambda):=g^{\lambda}(q)/f^{\lambda}.\]
Note that by the hook-content formula we have
\[L(\lambda^{\downarrow})/L(\lambda)=\frac{\prod_{u\in\lambda^{\downarrow}}\frac{q+c(u)}{i!}}{\prod_{u\in\lambda}\frac{q+c(u)}{(i-1)!}}
=\frac{q-\lambda_1'}{i}\]
where $\lambda_1'$ is the length of the first column of $\lambda$.

Summarizing, we have
\begin{equation}
\label{eqn:summ1}
E(\psi_i)^2\leq \frac{1}{q}\sum_{\lambda}\nu_{i,q}(\lambda^{\downarrow})
\frac{q-\lambda_1'}{i}=\frac{1}{qi}\left(q-\gamma_i\right),
\end{equation}
where $\gamma_i$ denotes the expectation of $(\lambda_1^{(i)})'$, i.e., the
expected length of the first column of a random shape with $i$ boxes, drawn
under the Plancherel-{\tt RSK} measure.

Notice also that since $\psi_i$ is an indicator random variable, we have
\begin{equation}
\label{eqn:summ2}
E(\psi_i^2)=E(\psi_i)=\gamma_i-\gamma_{i-1}.
\end{equation}
Therefore, combining (\ref{eqn:summ1}) and (\ref{eqn:summ2}) we obtain, by the
Cauchy-Schwarz inequality, the following difference inequality:
\begin{equation}
\label{eqn:crucial1}
\gamma_i-\gamma_{i-1}\leq \sqrt{\frac{1}{qi}}\sqrt{q-\gamma_i}.
\end{equation}
We claim that $\gamma_i\leq \beta(i)2\sqrt n$.

To prove this, note the following facts about $\gamma_i$:
\begin{itemize}
\item[(a)] $\gamma_{i+1}\geq \gamma_i$; and
\item[(b)] $\gamma_i\leq q$.
\end{itemize}
Now define a linear interpolation: for $t\in [i/q,(i+1)/q]$, set
\begin{equation}
\label{eqn:crucial2}
\beta_t=\frac{\gamma_i}{q}+q\left(t-\frac{i}{q}\right)
\left(\frac{\gamma_{i+1}}{q}-\frac{\gamma_{i}}{q}\right)\,.
\end{equation}
Note that for such $t$,
$\sqrt{1-\gamma_i/q}\leq \sqrt{1-\beta_t}$. Therefore, combining
equations (\ref{eqn:crucial1}) and (\ref{eqn:crucial2}) we obtain 
\begin{equation}
\frac{d}{dt} \beta_t=\gamma_{i+1}-\gamma_i\leq \frac{1}{\sqrt{qt}}\sqrt{1-\beta_t}\,,
\quad \beta_{1/q}=1/q\,. 
\end{equation}
Since $\gamma_i,\gamma_{i+1}\leq q$ we have $\beta_t\leq 1$, hence
the above differential inequality
is equivalent to
$$ -\frac{d}{dt}\sqrt{1-\beta_t}
\leq 1/(2\sqrt{qt})\,.$$
Hence it follows that
$$\sqrt{1-\beta_t}\geq \sqrt{1-1/q}-\sqrt{t/q}+1/q\geq
 \sqrt{1-1/q}-\sqrt{t/q}
\,.$$
That is,
$$\beta_t\leq 2\sqrt{t/q}-(t-1)/q\,.$$
Now we care about $t=n/q=\sqrt{n}/k$. We always have
$\beta_{n/q}\leq 1$ (trivially), but for $k>1$ we have the
better inequality $\beta_{n/q}\leq
2/k-1/k^2+1/q$. Therefore,
$$\gamma_n=q\beta_{n/q}\leq q\,,$$
for $k\leq 1$ and
$$\gamma_n=q\beta_{n/q}\leq (2-1/k)\sqrt{n}+1\,,$$
for $k>1$. The result then follows.\qed

\section*{Acknowledgments}
HT was supported by an NSERC Discovery Grant.
AY was supported by NSF grant DMS 0601010 and a U.~Minnesota DTC
grant during the Spring 2007; he also utilized the resources of
the Fields Institute, and of Algorithmics Incorporated, in Toronto, 
while a visitor. We would like
to thank Ofer Zeitouni for allowing us to include his proof for the
$\alpha>\alpha_{\rm critical}=\frac{1}{2}$ case of Theorem~\ref{thm:main} and for his
extensive help during this project.
We also thank Alexander Barvinok, Nantel Bergeron, Alexei Borodin,
Sergey Fomin, Christian Houdr\'{e}, Nicolas Lanchier, Igor Pak, Eric Rains,
Mark Shimozono,
Richard Stanley, Dennis Stanton, Craig Tracy and Alexander Woo 
for helpful correspondence.



\end{document}